\documentclass[11pt,twoside, epsf]{article}

\usepackage{color,graphicx,times}

\makeatother

\newcommand{\Across}{\raisebox{-0.25\height}{\includegraphics[width=0.5cm]{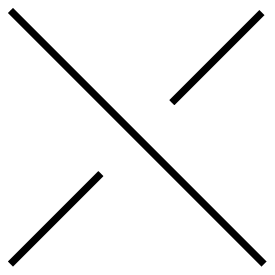}}}

\newcommand{\Asmooth}{\raisebox{-0.25\height}{\includegraphics[width=0.5cm]{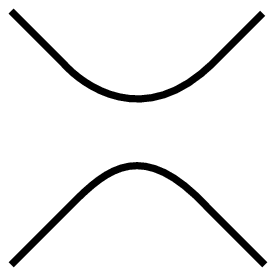}}}
\newcommand{\Bsmooth}{\raisebox{-0.25\height}{\includegraphics[width=0.5cm]{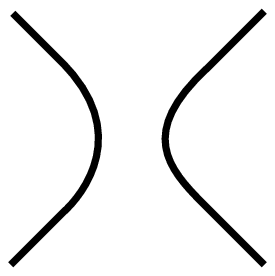}}}
\newcommand{\Rcurl}{\raisebox{-0.25\height}{\includegraphics[width=0.5cm]{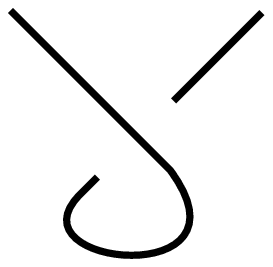}}}
\newcommand{\Lcurl}{\raisebox{-0.25\height}{\includegraphics[width=0.5cm]{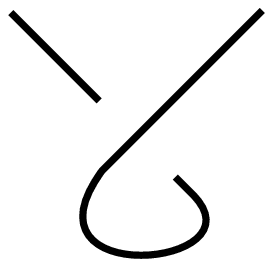}}}
\newcommand{\Arc}{\raisebox{-0.25\height}{\includegraphics[width=0.5cm]{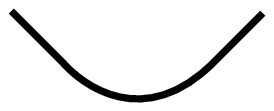}}}

\date{}

\title{Virtual Knot Cobordism}

{\author{Louis Hirsch Kauffman\\ Department of Mathematics, Statistics \\ and Computer Science (m/c 249)    \\ 851 South Morgan Street   \\ University of Illinois at Chicago\\
Chicago, Illinois 60607-7045\\ $<$kauffman@uic.edu$>$\\
}

\begin{document}

\maketitle

\abstract{ This paper defines a theory of cobordism for virtual knots and studies this theory for standard and rotational virtuals. Non-trivial examples of virtual slice knots are given and determinations of the four-ball genus of positive virtual knots are explained in relation to \cite{DKK}. Problems related to band-passing are explained and a theory of isotopy of virtual surfaces is formulated in terms of a generalization of the Yoshikawa moves.}
\bigbreak

\bigbreak

\bigbreak

\section{Introduction}
This paper defines a theory of cobordism for virtual knots. It is organized as follows.
In Section 2 we include a description of basics in virtual knot theory and the problems that arise from it. 
This section includes different interpretations of virtual knot theory including rotational virtual knot theory
(where the detour moves are constrained to be regular homotopies in the plane), problems related to the bracket and Jones polynomial for the theory, a discussion of parity, the odd writhe and  a description of the parity bracket polynomial formulated by Manturov \cite{MP}. All of this background material is used in the remainder of the paper where we apply these ideas and techniques to virtual knot cobordism. Section 3 gives the definitions for cobordism of virtual knots  and the definition of virtual Seifert surfaces. We define the four-ball genus $g_{4}(K)$ for virtual knots and links, and show that every virtual link $K$ bounds a virtual surface that is a natural generalization of the Seifert surface for a classical link. We state a new result about the four-ball genus for positive virtual knots. This result is proved in our paper \cite{DKK}. We give many properties of a key example, the virtual stevedore's knot. Then we discuss cobordism for rotational virtual knots and show that a virtual rotational link bounds a virtual surface if and only if it has an even number of viirtual crossings. Section 3.4 is an exploration of band-passing for virtual knots, and the possibility  to generalize classical results about the Arf invariant.
In Section 4 we develop a theory of virtual surfaces in four space based on a generalization of the Yoshikawa moves, and we give examples of virtual $2$-spheres in  four-space that are related to the virtual stevedore's knot. \bigbreak

This paper initiates a study of virtual knot cobordism at the level of virtual knot theory, a theory based on looking at knot diagrams without demanding planarity of the diagram.
Each crossing is locally in the plane, but the connections among the crossings can leave the plane and so arcs appear to cross over one another to make these connections. These crossings are called {\it virtual}, and the theory is expressed in terms of diagrams that contain both classical and virtual crossings. There are a number of different ways to interpret virtual knots and links topologically. We  describe these ways in Section 2. One can also develop more combinatorial approaches to virtual knot theory by giving up even more structure. One can use Gauss codes or Gauss diagrams to represent  virtual knots, and release certain structures related to the codes to make combinatorial theories that inform the virtual knot. Such work has been initiated by Turaev \cite{TurCob} and significantly by 
Manturov \cite{MBook,PCFree,CobFree}. We intend to connect the work of this paper with the work of Manturov in a joint paper under preparation. We should also mention the following papers on virtual knot theory that are useful background, but not cited directly in the present paper 
\cite{Kup,VKT,SA,DKT,DVK,Affine,ExtBr,Arrow,Diagram,Intro}.
\bigbreak

It should be remarked that the work \cite{PCFree}  can be used to prove
that many virtual knots are not concordant to any slice classical knot. This may not be immediately apparent to many readers, as the paper is focused on the cobordism of free knots. Results about 
free knots (undecorated  Gauss diagrams taken up Reidemeister move equivalence) are often applicable to standard virtual knots by simply forgetting some of the structure.  A crucial question about 
our formulation of virtual knot cobordism is: If two classical knots are concordant as virtual ones, are they concordant in the usual sense? At this writing, we do not know the answer to this question.
\bigbreak

It is worth mentioning areas of low dimensional topology and classical knot theory that are related to and informed by virtual knot theory. Just as classical combinatorial knot theory can be studied up to regular isotopy, it is useful to study virtual knots up to {\it rotational isotopy} (see Section 2.3 of the present paper). Then we have:
\begin{enumerate}
\item Every quantum invariant of classical knots extends to an invariant of rotational virtual knots.
\item Many quantum invariants, including the Jones polynomial, extend naturally to invariants of 
standard virtual knots. In the case of the Jones polynomial, it is inherent in its structure that it extends in a number of ways to invariants of virtual knots and rotational virtual knots. This means that the Jones
polynomial and these extensions becomes one of the main motivations for studying virtual knot theory.
The virtual knot theory becomes a way to ask new questions about the nature of the Jones polynomial.
\item Many problems related to categorification and link homology can be formulated in this domain.
\item The generalization of the classical braid group to the virtual braid group fits naturally into relationships of knot theory with quantum link invariants. In fact, in the context of the algebraic Yang-Baxter equation, the virtual braid group arises naturally for algebraic reasons (see \cite{CVBraid}).
\item Finally, the virtual knot theory has a number of variations. Standard virtual knot theory is equivalent to stabilized knot theory in thickened surfaces. When we add one of the forbidden moves, the scene changes and we have {\it welded knot theory} which is related to the Dahm braid group of circles in three-dimensional space and to embeddings of tori in four dimensional space. Then if we let go of structure we get flat virtual knot theory which is equivalent to stabilized immersions of curves in surfaces, and, letting go of even more structure, we have free virtual knot theory which is the study of Gauss diagrams without any orientation up to the Reidemeister moves. All of these subjects are interrelated and the combinatorial approach, combined with geometric topology makes a rich mixture of problems and ideas.
\end{enumerate}
\bigbreak

In this paper we formulate a theory of cobordism and concordance of virtual knots and links that is a generalization of the theory of cobordism of classical knots. This includes a virtual analog of ambient isotopy for surfaces in four-space and corresponding invariants of these isotopies. We give examples 
showing the viability of the generalizations and the problems that ensue.

 \section{Virtual Knot Theory}
Knot theory
studies the embeddings of curves in three-dimensional space.  Virtual knot theory studies the  embeddings of curves in thickened surfaces of arbitrary
genus, up to the addition and removal of empty handles from the surface. Virtual knots have a special diagrammatic theory, described below. Classical knot
theory embeds in virtual knot theory.
\bigbreak  

\begin{figure}[htb]
     \begin{center}
     \begin{tabular}{c}
     \includegraphics[width=8cm]{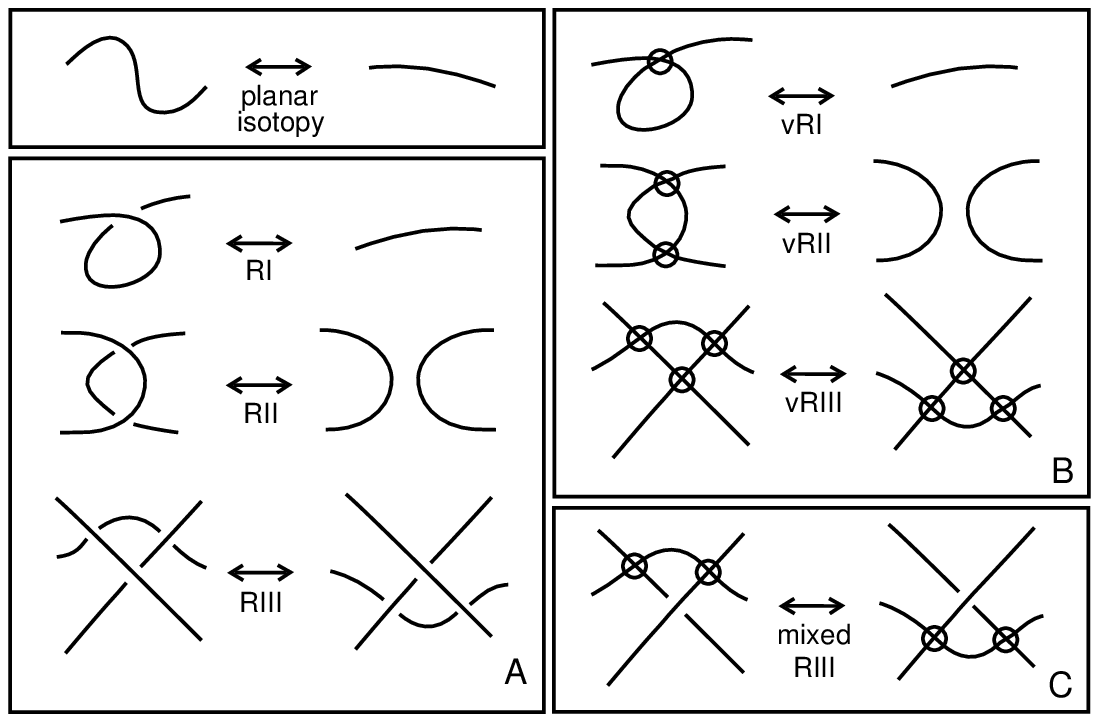}
     \end{tabular}
     \caption{\bf Moves}
     \label{Figure 1}
\end{center}
\end{figure}

In the diagrammatic theory of virtual knots one adds 
a {\em virtual crossing} (see Figure~\ref{Figure 1}) that is neither an over-crossing
nor an under-crossing.  A virtual crossing is represented by two crossing segments with a small circle
placed around the crossing point. 
\bigbreak

\begin{figure}[htb]
     \begin{center}
     \begin{tabular}{c}
     \includegraphics[width=6cm]{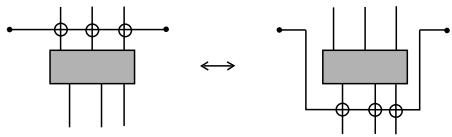}
     \end{tabular}
     \caption{\bf Detour Move}
     \label{Figure 2}
\end{center}
\end{figure}

Moves on virtual diagrams generalize the Reidemeister moves for classical knot and link
diagrams.  See  Figure~\ref{Figure 1}.  Classical crossings interact with
one another according to the usual Reidemeister moves, while virtual crossings are artifacts of the structure in the plane. 
Adding the global detour move to the Reidemeister moves completes the description of moves on virtual diagrams. In  Figure~\ref{Figure 1} we illustrate a set of local
moves involving virtual crossings. The global detour move is
	a consequence of  moves (B) and (C) in  Figure~\ref{Figure 1}. The detour move is illustrated in  Figure~\ref{Figure 2}.  Virtual knot and link diagrams that can be connected by a finite 
sequence of these moves are said to be {\it equivalent} or {\it virtually isotopic}.
\bigbreak

\begin{figure}[htb]
     \begin{center}
     \begin{tabular}{c}
     \includegraphics[width=6cm]{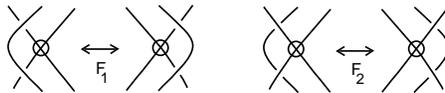}
     \end{tabular}
     \caption{\bf Forbidden Moves}
     \label{Figure 3}
\end{center}
\end{figure}

Another way to understand virtual diagrams is to regard them as representatives for oriented Gauss codes \cite{GPV}, \cite{VKT,SVKT} 
(Gauss diagrams). Such codes do not always have planar realizations. An attempt to embed such a code in the plane
leads to the production of the virtual crossings. The detour move makes the particular choice of virtual crossings 
irrelevant. {\it Virtual isotopy is the same as the equivalence relation generated on the collection
of oriented Gauss codes by abstract Reidemeister moves on these codes.}  
\bigbreak

 Figure~\ref{Figure 3} illustrates the two {\it forbidden moves}. Neither of these follows from Reidmeister moves plus detour move, and 
indeed it is not hard to construct examples of virtual knots that are non-trivial, but will become unknotted on the application of 
one or both of the forbidden moves. The forbidden moves change the structure of the Gauss code and, if desired, must be 
considered separately from the virtual knot theory proper. 
\bigbreak

\subsection{Interpretation of Virtuals Links as Stable Classes of Links in  Thickened Surfaces}
There is a useful topological interpretation \cite{VKT,DVK} for this virtual theory in terms of embeddings of links
in thickened surfaces.  Regard each 
virtual crossing as a shorthand for a detour of one of the arcs in the crossing through a 1-handle
that has been attached to the 2-sphere of the original diagram.  
By interpreting each virtual crossing in this way, we
obtain an embedding of a collection of circles into a thickened surface  $S_{g} \times R$ where $g$ is the 
number of virtual crossings in the original diagram $L$, $S_{g}$ is a compact oriented surface of genus $g$
and $R$ denotes the real line.  We say that two such surface embeddings are
{\em stably equivalent} if one can be obtained from another by isotopy in the thickened surfaces, 
homeomorphisms of the surfaces and the addition or subtraction of empty handles (i.e. the knot does not go through the handle).

\begin{figure}[htb]
     \begin{center}
     \begin{tabular}{c}
     \includegraphics[width=8cm]{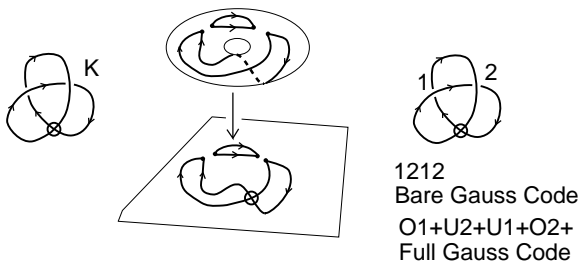}
     \end{tabular}
     \caption{\bf Surface Representation}
     \label{Figure 4}
\end{center}
\end{figure}

\noindent We have the
\smallbreak
\noindent
{\bf Theorem 1 \cite{VKT,DKT,DVK,Carter}.} {\em Two virtual link diagrams are isotopic if and only if their corresponding 
surface embeddings are stably equivalent.}  
\smallbreak
\noindent
\bigbreak  

\noindent  The reader will find more information about this correspondence \cite{VKT,DKT} in other papers by the author and in the literature of virtual knot theory.
\bigbreak
 
 \subsection{Review of the Bracket Polynomial for Virtual Knots}

In this section we recall how the bracket state summation model \cite{KaB} for the Jones polynomial  is defined for virtual knots
and links.  
\bigbreak

The bracket polynomial \cite{KaB} model for the Jones polynomial \cite{JO,JO1,JO2,Witten} is usually described by the expansion
\begin{equation}
\langle \Across \rangle=A \langle \Asmooth \rangle + A^{-1}\langle
\Bsmooth \rangle \label{kabr}
\end{equation}

and we have

\begin{equation}
\langle K \, \bigcirc \rangle=(-A^{2} -A^{-2}) \langle K \rangle \label{kabr1}
\end{equation}

\begin{equation}
\langle \Rcurl \rangle=(-A^{3}) \langle \Arc \rangle \label{kabr2}
\end{equation}

\begin{equation}
\langle \Lcurl \rangle=(-A^{-3}) \langle \Arc \rangle \label{kabr3}
\end{equation}
\bigbreak

We call a diagram in the plane 
{\em purely virtual} if the only crossings in the diagram are virtual crossings. Each purely virtual diagram is equivalent by the
virtual moves to a disjoint collection of circles in the plane.
\bigbreak

A state $S$ of a link diagram $K$ is obtained by
choosing a smoothing for each crossing in the diagram and labelling that smoothing with either $A$ or $A^{-1}$
according to the convention indicated in the bracket expansion above.  Then, given
a state $S$, one has the evaluation $<K|S>$ equal to the product of the labels at the smoothings, and one has the 
evaluation $||S||$ equal to the number of loops in the state (the smoothings produce purely virtual diagrams).  One then has
the formula
$$<K> = \Sigma_{S}<K|S>d^{||S||-1}$$
where the summation runs over the states $S$ of the diagram $K$, and $d = -A^{2} - A^{-2}.$
This state summation is invariant under all classical and virtual moves except the first Reidemeister move.
The bracket polynomial is normalized to an
invariant $f_{K}(A)$ of all the moves by the formula  $f_{K}(A) = (-A^{3})^{-w(K)}<K>$ where $w(K)$ is the
writhe of the (now) oriented diagram $K$. The writhe is the sum of the orientation signs ($\pm 1)$ of the 
crossings of the diagram. The Jones polynomial, $V_{K}(t)$ is given in terms of this model by the formula
$$V_{K}(t) = f_{K}(t^{-1/4}).$$
\noindent This definition is a direct generalization to the virtual category of the  
state sum model for the original Jones polynomial. It is straightforward to verify the invariances stated above.
In this way one has the Jones polynomial for virtual knots and links.
\bigbreak

\noindent We have \cite{DVK} the  
\smallbreak
\noindent
{\bf Theorem.} {\em To each non-trivial
classical knot diagram of one component $K$ there is a corresponding  non-trivial virtual knot diagram $Virt(K)$ with unit
Jones polynomial.} 
\bigbreak

\begin{figure}[htb]
     \begin{center}
     \begin{tabular}{c}
     \includegraphics[width=7cm]{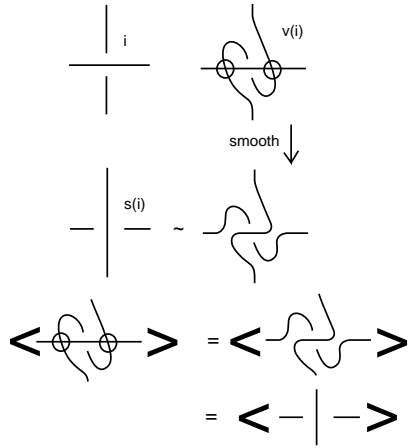}
     \end{tabular}
     \caption{\bf Virtualizing a Crossing and Crossing Switches}
     \label{Figure 5}
\end{center}
\end{figure}

\begin{figure}[htb]
     \begin{center}
     \begin{tabular}{c}
     \includegraphics[width=5cm]{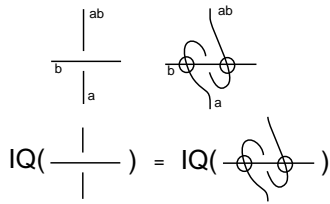}
     \end{tabular}
     \caption{\bf Quandle Invariance Under Virtualization}
     \label{Figure 6}
\end{center}
\end{figure}

The main ideas behind this Theorem are indicated in  Figure~\ref{Figure 5} and  Figure~\ref{Figure 6}.  In  Figure~\ref{Figure 5} we indicate the 
virtualization operation that replaces a classical crossing by using two virtual crossings and changing the implicit orientation of the classical crossing. We also show how the bracket polynomial sees this operation as though the crossing had been switched in the classical knot. Thus, if we virtualize as set
of classical crossings whose switching will unknot the knot, then the virtualized knot will have unit 
Jones polynomial. On the other hand, the virtualization is invisible to the quandle, as shown in  
Figure~\ref{Figure 6}. This implies (by properties of the quandle) that virtual knots obtained in this way from classical non-trivial knots will themselves be non-trivial.
\bigbreak

It is an open problem whether there are classical knots (actually knotted) having unit Jones polynomial. (There are linked links whose linkedness is unseen
\cite{EKT} by the Jones polynomial.) If there exists a classical knot with unit Jones polynomial, then one of the knots $Virt(K)$ produced by this Theorem
may be isotopic to  a classical knot.  Such examples are guaranteed to be non-trivial, but they are usually also not classical. 
We do not know at this writing whether all such virtualizations of non-trivial classical knots, yielding virtual knots with unit Jones polynomial, are 
non classical.  This has led to an investigation of new invariants for virtual knots. 
\bigbreak

\begin{figure}
     \begin{center}
     \begin{tabular}{c}
     \includegraphics[width=6cm]{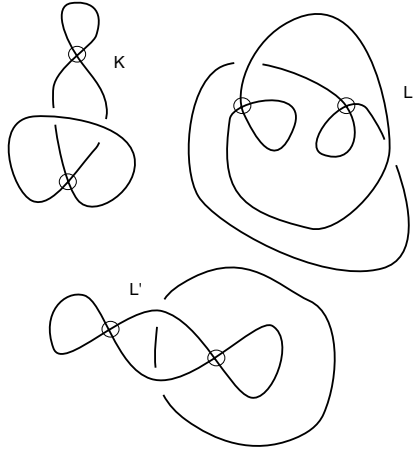}
     \end{tabular}
     \caption{\bf A rotational virtual knot and two rotational virtual links.}
     \label{rotateknot}
\end{center}
\end{figure}

%\begin{figure}
%     \begin{center}
%     \begin{tabular}{c}
%     \includegraphics[width=5cm]{RotLink.eps}
 %    \end{tabular}
 %    \caption{A rotational virtual knot and a rotational virtual link.}
 %    \label{rotateknot}
%\end{center}
%\end{figure}

\subsection{Rotational Virtual Knot Theory}
\begin{figure}
     \begin{center}
     \begin{tabular}{c}
     \includegraphics[width=6cm]{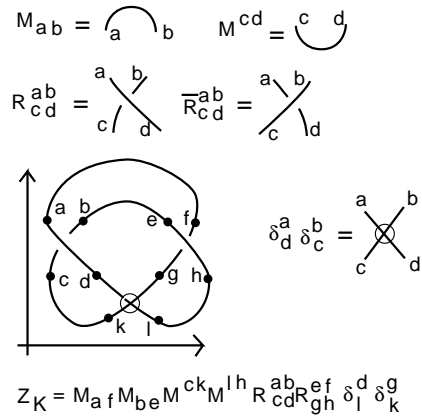}
     \end{tabular}
     \caption{\bf Quantum Link Invariants}
     \label{quantum}
\end{center}
\end{figure}

Rotational virtual knot theory introduced in~\cite{VKT} is virtual
knot theory without the first virtual move (thus one does not allow
the addition or deletion of a virtual curl. All diagrams are represented in the surface of plane so that we can distinguish clockwise from counterclockwise rotations.. 
See Figure~\ref{rotateknot} for an
example of a rotational virtual knot.
The rotational  version of virtual knot theory is significant because{\it  all quantum link invariants
originally defined for classical links extend to rotational virtual
knot theory.} We give a quick illustration of this point in Figure~\ref{quantum} where we indicate how a quantum link invariant depends on matrices or operators assigned to each crossing and each maximum and minimum of the diagram. One extends this to virtual crossings by using a crossed identity operator ( a transposition) at the vritual crossings. The possibly non-trivial maxima and minima create the need to use the rotational context. This theory has been explored~\cite{VKT,CVBraid} and deserves further exploration. 
\bigbreak

\begin{figure}
     \begin{center}
     \begin{tabular}{c}
     \includegraphics[width=8cm]{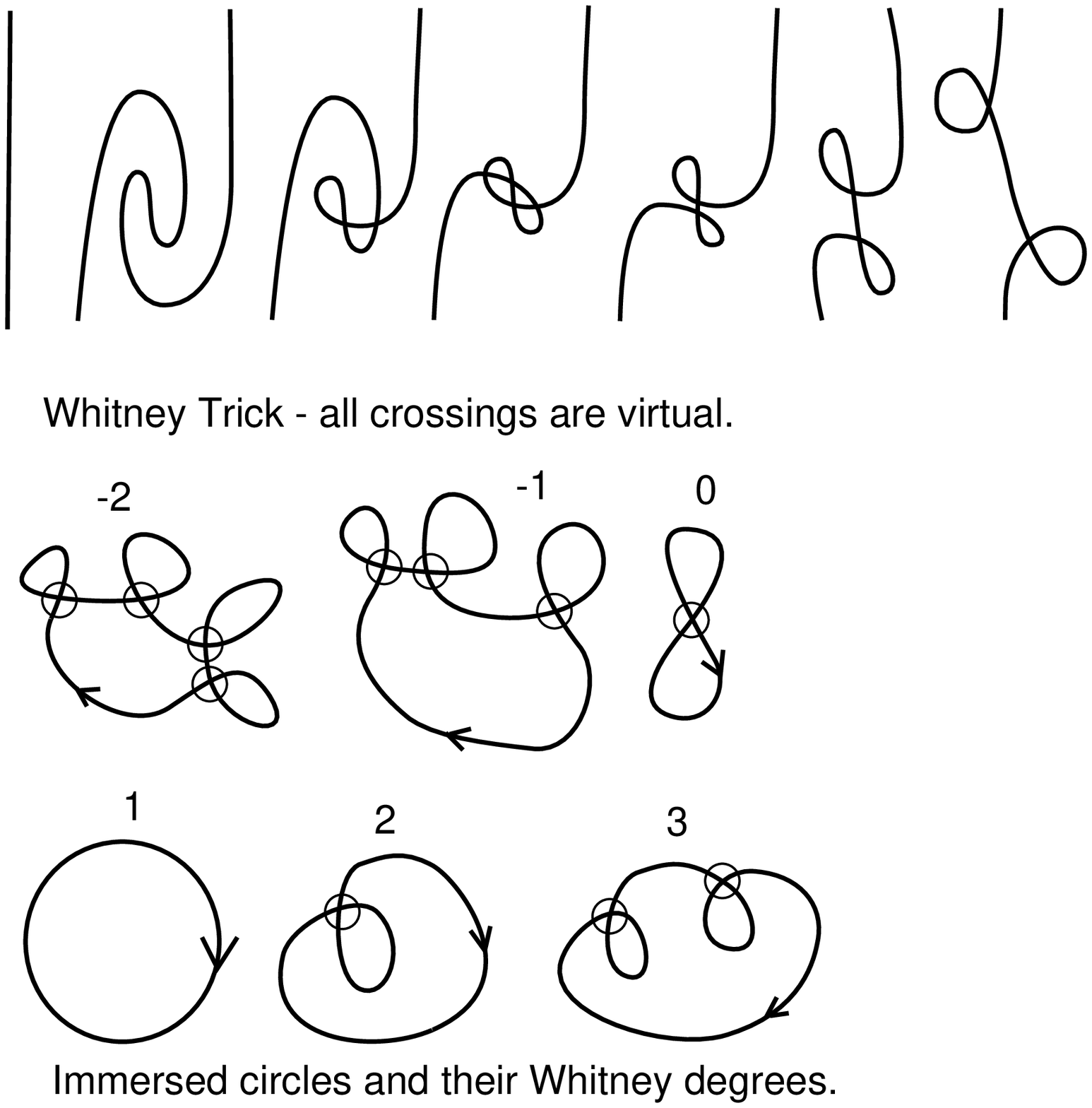}
     \end{tabular}
     \caption{\bf Whitney Trick and Whitney Degrees}
     \label{whitney}
\end{center}
\end{figure}

\begin{figure}
     \begin{center}
     \begin{tabular}{c}
     \includegraphics[width=6cm]{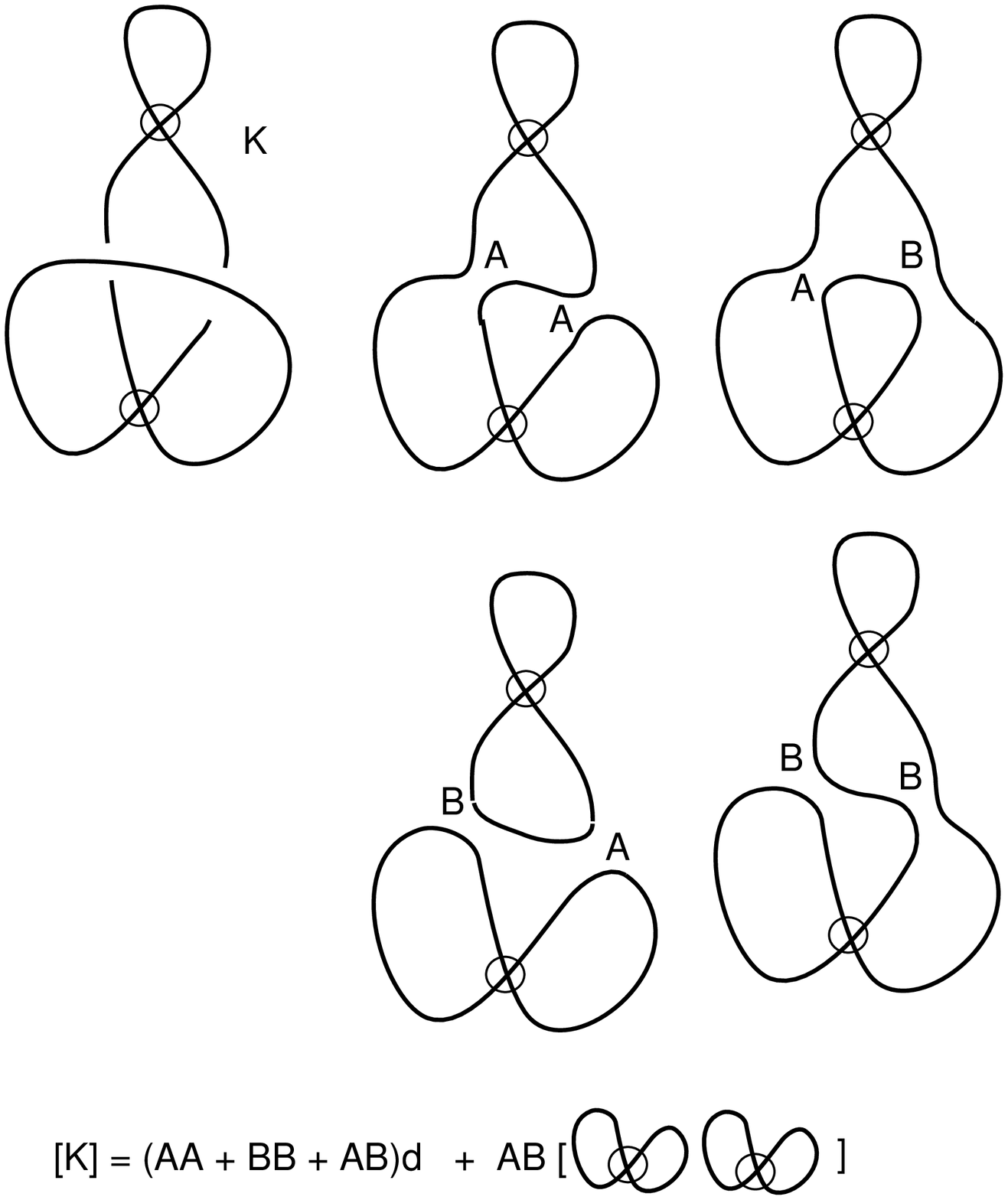}
     \end{tabular}
     \caption{\bf Bracket Expansion of a Rotational Virtual Knot}
     \label{rotbracket}
\end{center}
\end{figure}

We  formulate \cite{VKT} a version of the bracket polynomial for rotational
virtuals that assigns variables according to the Whitney degree of state curves. 
For rotational virtuals we extend the bracket just as we did for virtual knots and links except that the state curves are now disjoint unions in the plane of curves that have only virtual self-intersections, and are taken up to {\it regular homotopy} in the plane (we can take regular homotopy of curves to mean that the flat virtual versions of the second and third Reidemeister moves are allowed and that regular homotopy is the equivalence relation generated by these moves and planar homeomorphisms). When we expand the bracket we obtain a state sum of the form
$$[K] = \Sigma_{S}<K|S>[S]$$
where the summation is over all states obtained by smoothing every crossing in the virtual diagram $K$
and $<K|S>$ is the product of the weights $A$ and $A^{-1}$ just as before. An empty loop with no virtual crossings (in its virtual equivalence class) will be evaluated as $d = -A^2 - A^{-2}.$
The symbol  $[S]$ is the {\it planar class} of the state $S.$ By the planar class of the state we mean its equivalence class up to virtual rotational equivalence. This means that each state loop is taken as a regular homotopy class. These individual classes are in $1$-$1$ correspondence with the integers, as shown in Figure~\ref{whitney} (via the Whitney trick and the winding degree of the plane curves), and can be handled by using combinatorial regular isotopy as in \cite{FKT}.  A configuration of loops (possibly nested) is equivalent to a disjoint union of adjacent loops. We can thus regard each virtual loop as a variable $d_{n}$ where $n$ is an integer and $d_{1} =d_{-1} = -A^2 - A^{-2}.$ Here we give an examples of a computation of 
$[K]$ for a rotational virtual knots in Figure~\ref{rotbracket}. The reader will note that in this example, even if we let $A = -1 = B$ and $d = -2$ the invariant is still non-trivial due to the appearance of the two loops with Whitney degree zero. Thus the example in Figure~\ref{rotbracket} also gives a non-trivial flat rotational virtual knot. We shall look at cobordism of rotational virtual knots later in the paper.
\bigbreak
 
\subsection{Parity, Odd Writhe and the Parity Bracket Polynomial}
Parity is an important theme in virtual knot theory and figures in many investigations of this subject.
In a virtual knot diagram there can be both even and odd crossings.  A crossing is {\it odd} if 
it flanks an odd number of symbols in the Gauss code of the diagram. A crossing is {\it even} if 
it flanks an even number of symbols in the Gauss code of the diagram. For example, in 
Figure~\ref{Figure 4} we illustrate the virtual knot $K$ with bare Gauss code $1212.$ Both crossings in the diagram $K$ are odd. In any classical knot diagram all crossings are even.
\bigbreak 

In \cite{SL} we introduced the  {\it odd writhe} $J(K)$ for any virtual diagram $K.$ $J(K)$  is the sum of the signs of the odd crossings. Classical diagrams have zero odd writhe. Thus if $J(K)$ is non-zero, then $K$ is not equivalent to any classical knot. For the mirror image $K^{*}$ of any diagram $K,$ we have the formula $J(K^{*}) = - J(K).$ Thus, when $J(K) \ne 0,$  we know that the knot $K$ is
not classical and not equivalent to its mirror image. Parity does all the work in this simple invariant.
For example, if $K$ is the virtual knot in Figure~\ref{Figure 4}, the we have $J(K) = 2.$ Thus $K,$ the simplest virtual knot, is non-classical and it is chiral (inequivalent to its mirror image.)
\bigbreak
 
\bigbreak
In this section we introduce the Manturov Parity Bracket  \cite{MP}.
This is a form of the bracket polynomial defined for virtual knots and for free knots (unlabeled Gauss diagrams taken up to abstract Reidemeister move equivalence)  that uses the parity of the crossings. To compute the parity bracket, we first make all the odd crossings into graphical vertices. Then we expand the resulting diagram on the remaining even crossings. The result is a sum of graphs with polynomial coefficients. 
\bigbreak

\begin{figure}
     \begin{center}
     \begin{tabular}{c}
     \includegraphics[width=8cm]{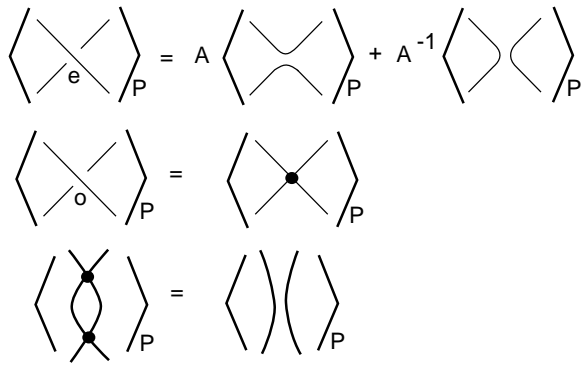}
     \end{tabular}
     \caption{\bf Parity Bracket Expansion}
     \label{pbracket}
\end{center}
\end{figure}

More precisely, let $K$ be a virtual knot diagram. Let $E(K)$ denote the result of making all the odd crossings in $K$ into graphical nodes as illustrated in  Figure~\ref{pbracket} .
Let $SE(K)$ denote the set of all bracket states of $E(K)$ obtained by smoothing each classical crossing in $E(K)$ in one of the two possible ways. Then we define the {\it parity bracket} 
$$<K>_{P} = (1/d)\Sigma_{S \in SE(K)} A^{i(S)} [S]$$ where $d=-A^2 - A^{-2}$, $i(S)$ denotes the 
product of $A$ or $A^{-1}$ from each smoothing site according to the conventions of  Figure~\ref{pbracket}, and $[S]$ denotes the reduced class of the virtual graph $S.$  The graphs are subject to a reduction move that eliminates bigons as in the second Reidemeister move on a knot diagram as shown in Figure~\ref{pbracket}.  Thus $[S]$ represents the unique minimal representative for the virtual graph $S$ under virtual graph isotopy coupled with the bigon reduction move. A graph that reduces to a circle (the circle is a graph for our purposes) is replaced by the value $d$ above. Thus $<K>_{P}$ is an element of a module  generated by reduced graphs with  coefficients Laurent polynomials in $A.$.
\bigbreak

\begin{figure}
     \begin{center}
     \begin{tabular}{c}
     \includegraphics[width=4cm]{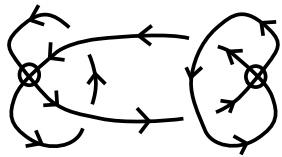}
     \end{tabular}
     \caption{\bf Kishino Diagram}
     \label{kishino}
\end{center}
\end{figure}

With the usual bracket polynomial variable $A$, the parity bracket is an invariant of standard virtual knots. With $A=\pm 1$ it is an invariant of flat virtual knots. Even more simply, with $A=1$ and taken modulo two, we have an invariant of flat knots with loop value zero. See Figure~\ref{kishino} for an illustration of the application of the parity bracket to the Kishino diagram illustrated there. The Kishino diagram is notorious for being hard to detect by the usual polynomial invariants such as the Jones polynomial. It is a perfect example of the power of the parity bracket. All the crossings of the Kishino diagram are odd. Thus there is exactly one term in the evaluation of the Kishino diagram by the parity bracket, and this term is the Kishino diagram itself, with its crossings made into graphical nodes. The resulting graph is irreducible and so the Kishino diagram becomes its own invariant. We conclude that this diagram will be found from any isotopic version of the Kishino diagram. This allows strong conclusions about many properties of the diagram. For example, it is easy to check that the least surface on which this diagram can be represented with the given planar cyclic orders at the nodes) is genus two. Thus we conclude that the least genus for a surface representation of the Kishino diagram as a flat knot or virtual knot is two.
\bigbreak

\begin{figure}
     \begin{center}
     \begin{tabular}{c}
     \includegraphics[width=4cm]{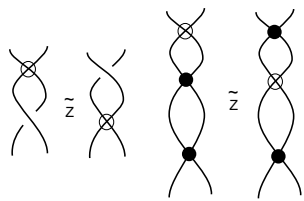}
     \end{tabular}
     \caption{\bf Z-Move and Graphical Z-Move}
     \label{Figure9}
\end{center}
\end{figure}

\begin{figure}
     \begin{center}
     \begin{tabular}{c}
     \includegraphics[width=6cm]{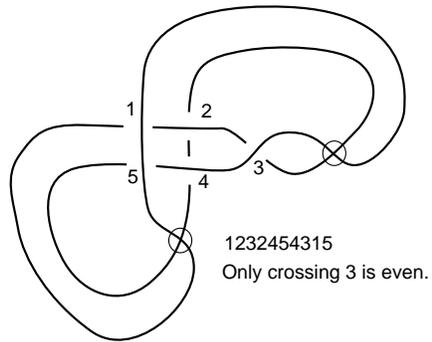}
     \end{tabular}
     \caption{\bf A Knot KS With Unit Jones Polynomial}
     \label{Figure10}
\end{center}
\end{figure}

\begin{figure}
     \begin{center}
     \begin{tabular}{c}
     \includegraphics[width=6cm]{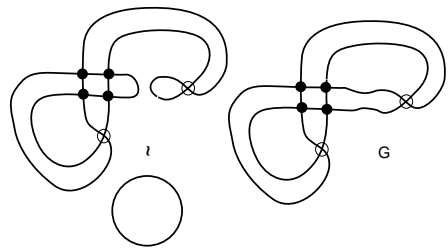}
     \end{tabular}
     \caption{\bf Parity Bracket States for the Knot KS}
     \label{Figure11}
\end{center}
\end{figure}

In Figure~\ref{Figure9} we illustrate the {\it Z-move}  and the {\it graphical Z-move}. Two virtual knots or links
that are related by a Z-move have the same standard bracket polynomial. This follows directly from our discussion in the previous section. We would like to analyze the structure of Z-moves using the parity
bracket. In order to do this we need a version of the parity bracket that is invariant under the Z-move.
In order to accomplish this, we need to add a corresponding Z-move in the graphical reduction process for the parity bracket. This extra graphical reduction is indicated in Figure~\ref{Figure9} where we show a graphical Z-move. The reader will note that graphs that are irreducible without the graphical Z-move can become reducible if we allow graphical Z-moves in the reduction process. For example, the graph associated with the Kishino knot is reducible under graphical Z-moves. However, there are examples of 
graphs that are not reducible under graphical Z-moves and Reidemeister two moves. An example of such a graph occurs in the parity bracket of the knot $KS$ shown in Figure~\ref{Figure10} and 
Figure~\ref{Figure11}. This knot has one even classical crossing and four odd crossings. One smoothing of the even crossing yields a state that reduces to a loop with no graphical nodes, while the other smoothing yields a state that is irreducible even when the Z-move is allowed. The upshot is that this knot KS is not Z-equivalent to any classical knot. Since one can verify that $KS$ has unit Jones polynomial, this example is a counterexample to a conjecture of 
Fenn, Kauffman and Maturov  \cite{FKM} that suggested that a knot with unit Jones polynomial should be Z-equivalent to a classical knot.
\bigbreak

\section{Virtual Knot Cobordism}
\begin{figure}
     \begin{center}
     \begin{tabular}{c}
     \includegraphics[width=8cm]{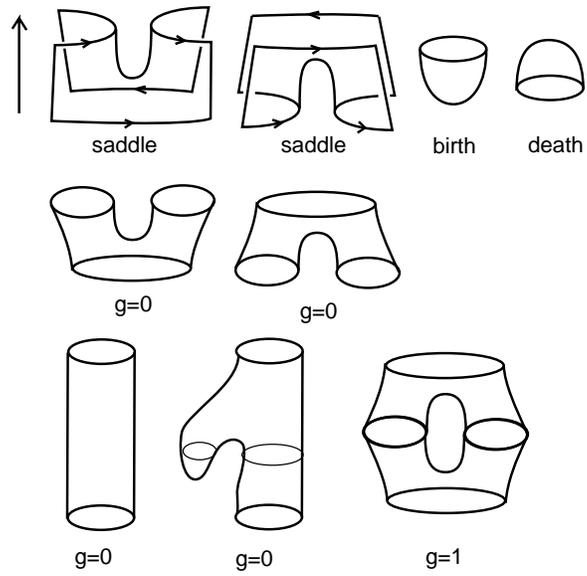}
     \end{tabular}
     \caption{\bf Saddles, Births and Deaths}
     \label{saddle}
\end{center}
\end{figure}

\noindent {\bf Definition.} Two oriented knots or links $K$ and $K'$ are {\it virtually cobordant}  if one may be obtained from the other by a sequence of virtual isotopies (Reidemeister moves plus detour moves) plus births, deaths and oriented saddle points, as illustrated in  Figure~\ref{saddle}. A {\it birth} is the introduction into the diagram of an isolated unknotted circle. A {\it death} is the removal from the diagram
of an isolated unknotted circle. A saddle point move results from bringing oppositely oriented arcs 
into proximity and resmoothing the resulting site to obtain two new oppositely oriented arcs. See the Figure for an illustration of the process. Figure\ref{saddle} also illustrates the {\it schema} of surfaces that are generated by  cobordism process. These are abstract surfaces with well defined genus in terms of the sequence of steps in the cobordism. In the Figure we illustrate two examples of genus zero, and one example of genus 1. We say that a cobordism has genus $g$ if its schema has that genus. Two knots are {\it cocordant} if there is a cobordism of genus zero  connecting them. A virtual knot is said to be a {\it slice} knot if it is virtually concordant to the unknot, or equivalently if it is virtually concordant to the empty knot (The unknot is concordant to the empty knot via one death.). As we shall see below, {\it every virtual knot or link is concordant to the unknot}. Another way to say this, is to say that there is a {\it virtual surface} (schema) whose boundary is the give virtual knot. The reader should note that when we speak of a virtual surface, we mean a surface schema that is generated by saddle moves, maxima and minima as describe above. 
\bigbreak 

\noindent{\bf Definition.} We define the {\it four-ball genus} $g_{4}(K)$ of a virtual knot or link $K$ to be the least genus among all virtual surfaces that bound $K.$ As we shall see below, there is a simple upper bound on the four-ball genus for any virtual knot or link and a definite result for the four-ball genus of positive virtual knots \cite{DKK}.
\bigbreak

\begin{figure}
     \begin{center}
     \begin{tabular}{c}
     \includegraphics[width=8cm]{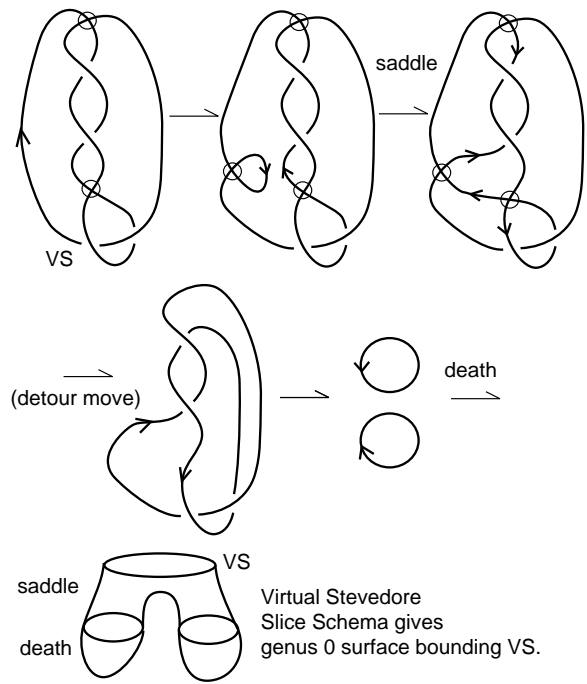}
     \end{tabular}
     \caption{\bf Virtual Stevedore is Slice}
     \label{vstevedore}
\end{center}
\end{figure}

In Figure~\ref{vstevedore} we illustrate the {\it virtual stevedore's knot, VS} and show that it is a slice knot in the sense of the above definition. This figure illustrates how the surface schema whose boundary in the virtual stevedore is evolved via the saddle point that produces two virtually unlinked curves that are isotopic to a pair of curves that can undergo deaths to produce the genus zero slicing surface. We will use this example to illustrate our theory of virtual knot cobordism, and the questions that we are investigating. Before looking at the virtual stevedore in this detail, we make a digression about spanning surfaces and the four-ball genus.
\bigbreak

\bigbreak

\subsection{Spanning Surfaces for Knots and Virtual Knots.}
\begin{figure}
     \begin{center}
     \begin{tabular}{c}
     \includegraphics[width=7cm]{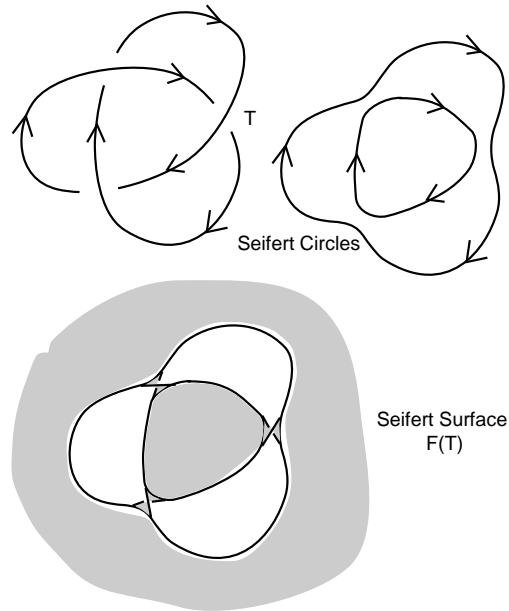}
     \end{tabular}
     \caption{\bf Classical Seifert Surface}
     \label{seifert}
\end{center}
\end{figure}

It is a well-known that every oriented classical knot or link bounds an embedded orientable surface in three-space. A representative surface of this kind can be obtained by the algorithm due to Seifert (See \cite{OK,FKT,KP}). We have illustrated Seifert's algorithm for a trefoil diagram in Figure~\ref{seifert}. The algorithm proceeds as follows: At each oriented crossing in a given diagram 
$K,$ smooth that crossing in the oriented manner (reconnecting the arcs locally so that the crossing disappears and the connections respect the orientation). The result of this operation is a collection of oriented simple closed curves in the plane, usually called the {\it Seifert circles}. To form the {\it Seifert surface} $F(K)$ for the diagram $K,$ attach disjoint discs to each of the Seifert circles, and connect these discs to one another by local half-twisted bands at the sites of the smoothing of the diagram. This process is indicated in the Figure~\ref{seifert}. In that figure we have not completed the illustration of the outer disc.
\bigbreak

It is important to observe that we can calculate the genus of the resulting surface quite easily from the combinatorics of the classical knot diagram $K.$ For purposes of simplicity, we shall assume that we are dealing with a knot diagram (one boundary component) and leave the case of links to the reader.
We then have the 
\bigbreak

\noindent{\bf Lemma.} Let $K$ be a classical knot diagram with $n$ crossings and $r$ Seifert circles.
then the genus of the Seifert Surface $F(K)$ is given by the formula
$$g(F(K)) =(1/2)( -r + n +1).$$
\bigbreak

\noindent {\bf Proof.} The surface $F(K),$ as described prior to the statement of the Lemma, retracts to a cell complex consisting of the projected graph of the knot diagram with two-cells attached to each cycle in the graph that corresponds to a Seifert circle. Thus we have that the Euler characteristic of this surface is given the the formula $$\chi(F(K)) = n - e + r$$ where $n,$ the number of crossings in the diagram, is the number of zero-cells, $e$ is the number of one-cells (edges) in the projected diagram (from node to node), and $r$ is the number of Seifert circles as these are in correspondence with the two-cells. However, we know that $4n = 2e$ since there are four edges locally incident to each crossing. Thus,
$$\chi(F(K)) =  - n + r.$$ Furthermore, we have that $\chi(F(K)) = 1 - 2g(F(K)),$ since  this surface has a single boundary component and is orientable. From this it follows that $1-2g(F(K)) = -n + r,$ and hence
$$g(F(K)) =(1/2)( -r + n +1).$$  This completes the proof. //
\bigbreak

\begin{figure}
     \begin{center}
     \begin{tabular}{c}
     \includegraphics[width=7cm]{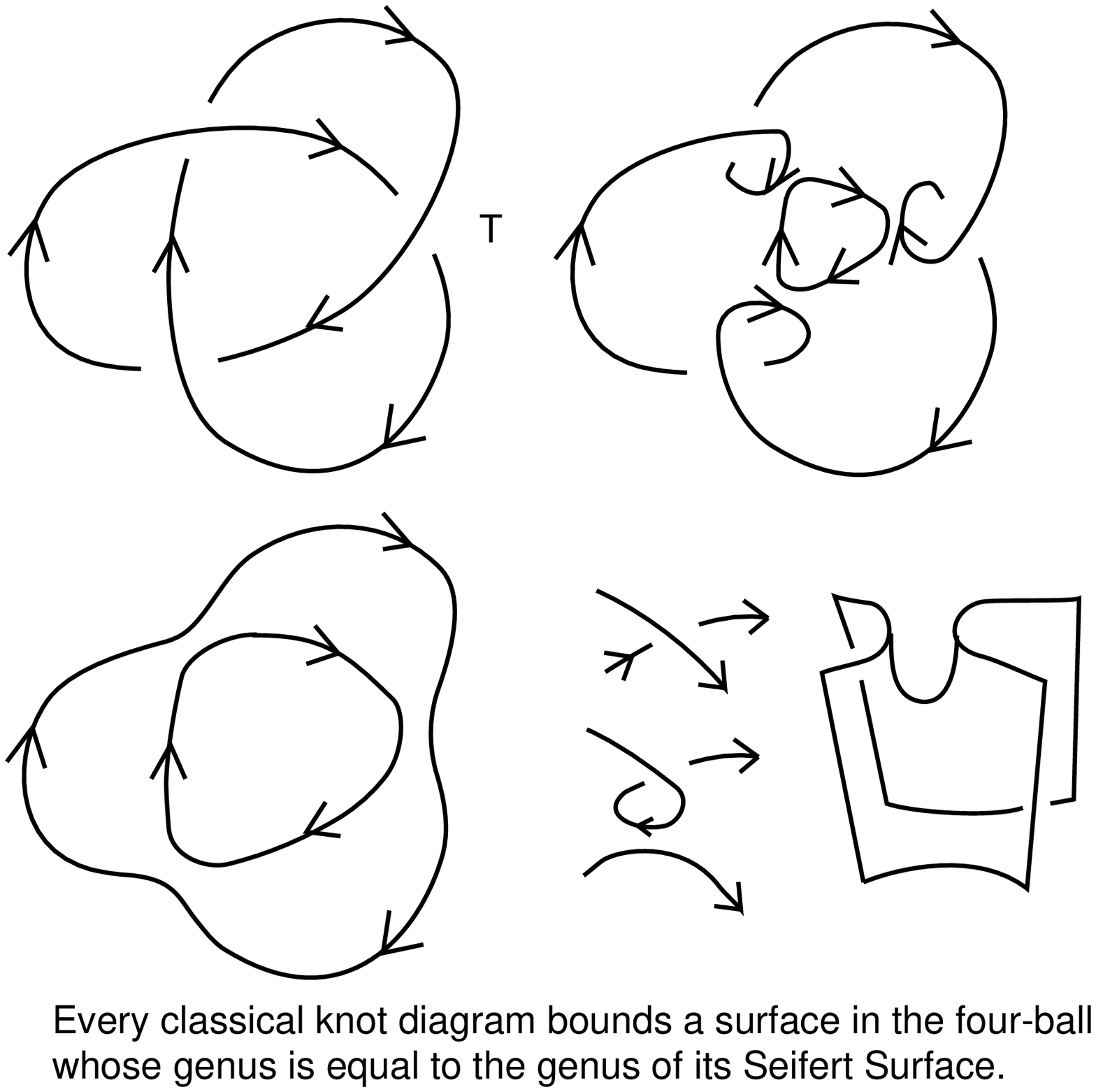}
     \end{tabular}
     \caption{\bf Classical Cobordism Surface}
     \label{classicalcob}
\end{center}
\end{figure}

We now observe that {\it for any classical knot $K,$ there is a surface bounding that knot in the four-ball that is homeomorphic to the Seifert surface}.  One can construct this surface by  pushing the Seifert
surface into the four-ball keeping it fixed along the boundary. We will give here a  different description of
this surface as indicated in Figure~\ref{classicalcob}. In that figure we {\it perform a saddle point transformation at every crossing of the diagram.} The result is a collection of unknotted and unlinked curves. By our interpretation of surfaces in the four-ball obtained by saddle moves and isotopies, we can then bound each of these curves by discs (via deaths of circles) and obtain a surface 
$S(K)$ embedded in the four-ball with boundary $K.$ As the reader can easily see, the curves produced by the saddle transformations are in one-to-one correspondence with the Seifert circles for
$K,$ and it easy to verity that $S(K)$ is homeomorphic with the Seifert surface $F(K).$
Thus we know that $g(S(K)) =(1/2)( -r + n +1).$ In fact the same argument that we used to analyze the genus of the Seifert surface applies directly to the construction of $S(K)$ via saddles and minima.
\bigbreak

\begin{figure}
     \begin{center}
     \begin{tabular}{c}
     \includegraphics[width=8cm]{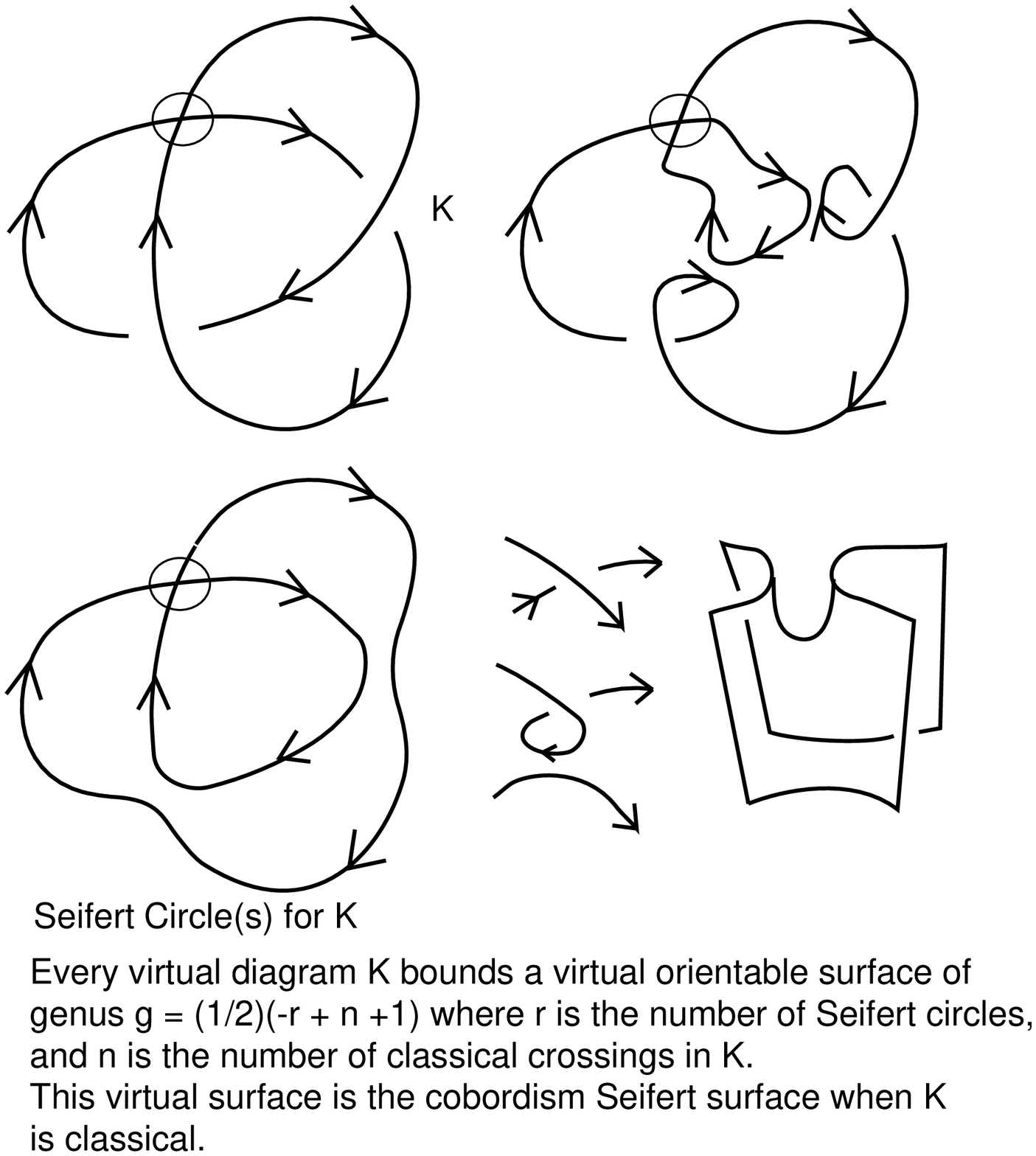}
     \end{tabular}
     \caption{\bf Virtual Cobordism Seifert Surface}
     \label{virtseifert}
\end{center}
\end{figure}

\begin{figure}
     \begin{center}
     \begin{tabular}{c}
     \includegraphics[width=7cm]{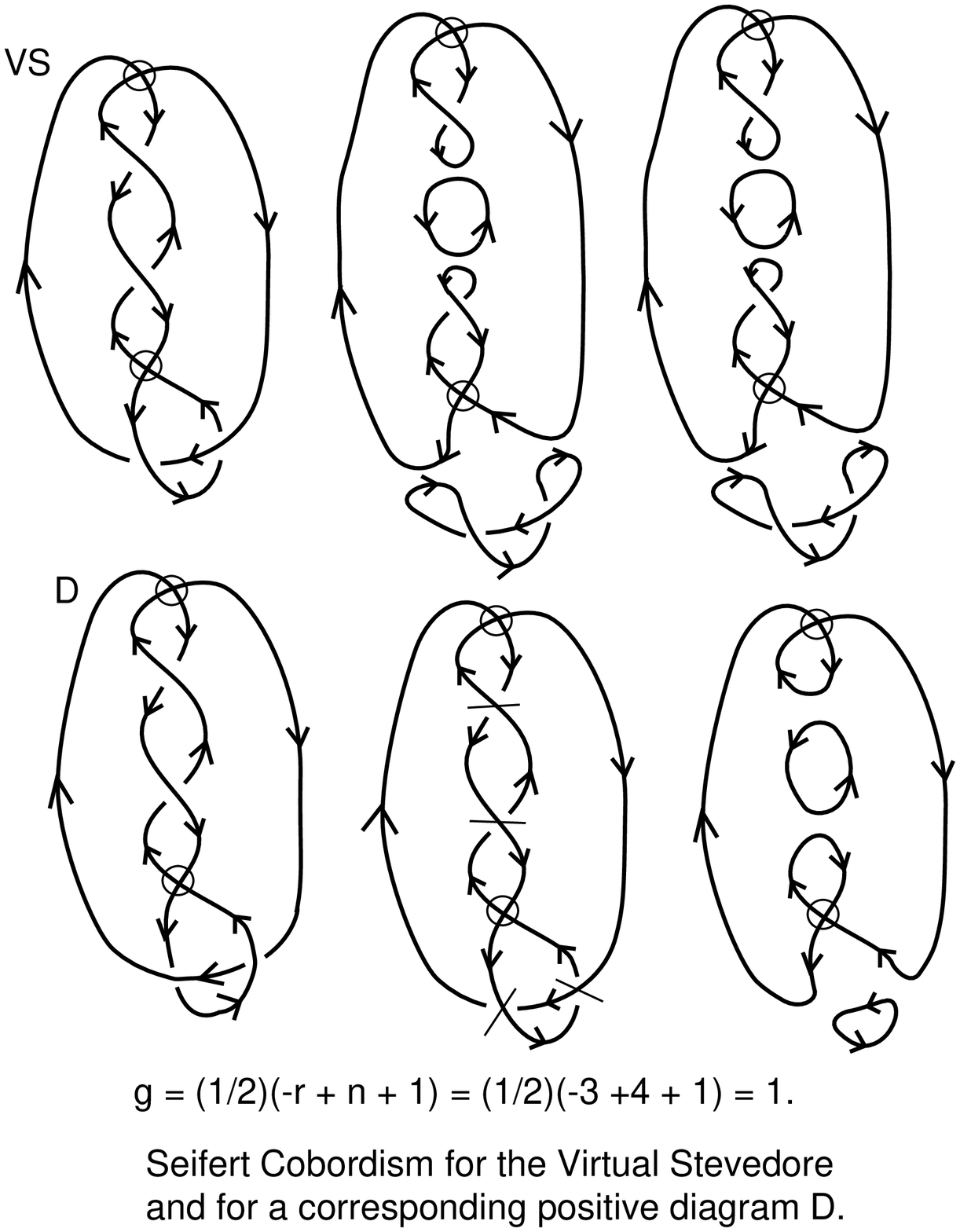}
     \end{tabular}
     \caption{\bf Virtual Stevedore Cobordism Seifert Surface}
     \label{vstevedoreseifert}
\end{center}
\end{figure}

Now the stage is set for generalizing the Seifert surface to a surface $S(K)$ for virtual knots $K.$
View Figure~\ref{virtseifert} and Figure~\ref{vstevedoreseifert}. In these figures we have performed a saddle transformation at each classical crossing of a virtual knot $K.$ The result is a collection of unknotted curves that are isotopic (by the first classical Reidemeister move) to curves with only virtual crossings. Once the first Reidemeister moves are performed, these curves are identical with the 
{\it virtual Seifert circles} obtained from the diagram $K$ by smoothing all of its classical crossings.
We can then Isotope these circles into a disjoint collection of circles (since they have no classical crossings) and cap them with discs in the four-ball. The result is a virtual surface $S(K)$ whose boundary is the given virtual knot $K.$ We will use the terminology {\it virtual surface in the four-ball}
for this surface schema. In the case of a virtual slice knot, we have that the knot bounds a virtual surface of genus zero. But with this construction we have proved the 
\bigbreak

\noindent{\bf Lemma.} Let $K$ be a virtual knot, then the virtual Seifert surface $S(K)$ constructed above
has genus given by the formula $$g(S(K)) = (1/2)(-r + n + 1)$$ where $r$ is the number of virtual Seifert
circles in the diagram $K$ and $n$ is the number of classical crossings in the diagram $K.$
\bigbreak

\noindent{\bf Proof.} The proof follows by the same argument that we already gave in the classical case.
Here the projected virtual diagram gives a four-regular graph $G$ (not necessarily planar) whose nodes are in one-to-one correspondence with the classical crossings of $K.$ The edges of $G$ are in one-to-one correspondence with the edges in the diagram that extend from one classical crossing to the next.
We regard $G$ as an abstract graph so the the virtual crossings disappear.
The argument then goes over verbatim in the sense that $G$ with two-cells attached to the virtual Seifert circles is a retract of the surface $S(K)$  constructed by cobordism. The counting argument for the genus is identical to the classical case. This completes the proof. //
\bigbreak

\noindent{\bf Remark.} For the virtual stevedore in Figure~\ref{vstevedoreseifert}
we have the interesting phenomenon that there is a much lower genus surface that can be produced by cobordism than the virtual Seifert surface. In that same figure we have illustrated a diagram $D$ with the same projected diagram as the virtual stevedore, but $D$ has all positive crossings. In this case we can prove \cite{DKK} that there is no virtual surface for this diagram  $D$ of four-ball genus less than $1.$ 
In fact, we have the following result which is proved in \cite{DKK}. This Theorem is a generalization of a corresponding result for classical knots due to Rasmussen \cite{Ras}.
\bigbreak

\noindent{\bf Theorem \cite{DKK}.} Let $K$ be a positive virtual knot (all classical crossings in $K$ are positive), then the four-ball genus $g_{4}(K)$ is given by the formula
$$g_{4}(K) = (1/2)(-r + n + 1) = g(S(K))$$ where $r$ is the number of virtual Seifert circles in the diagram
$K$ and $n$ is the number of classical crossings in this diagram. In other words, that virtual Seifert surface for $K$ represents its minimal four-ball genus.
\bigbreak

\begin{figure}
     \begin{center}
     \begin{tabular}{c}
     \includegraphics[width=7cm]{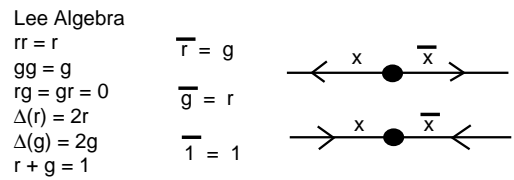}
     \end{tabular}
     \caption{\bf Lee Algebra Undergoes Involution at a Cut Locus}
     \label{leecut}
\end{center}
\end{figure}

\begin{figure}
     \begin{center}
     \begin{tabular}{c}
     \includegraphics[width=7cm]{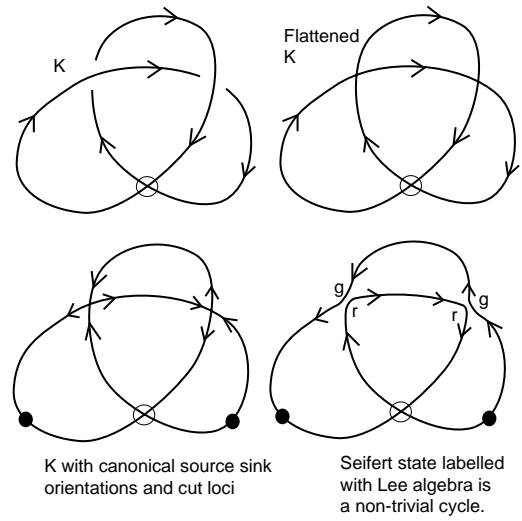}
     \end{tabular}
     \caption{\bf Lee Algebra Labels Seifert State}
     \label{leecycle}
\end{center}
\end{figure}

\begin{figure}
     \begin{center}
     \begin{tabular}{c}
     \includegraphics[width=7cm]{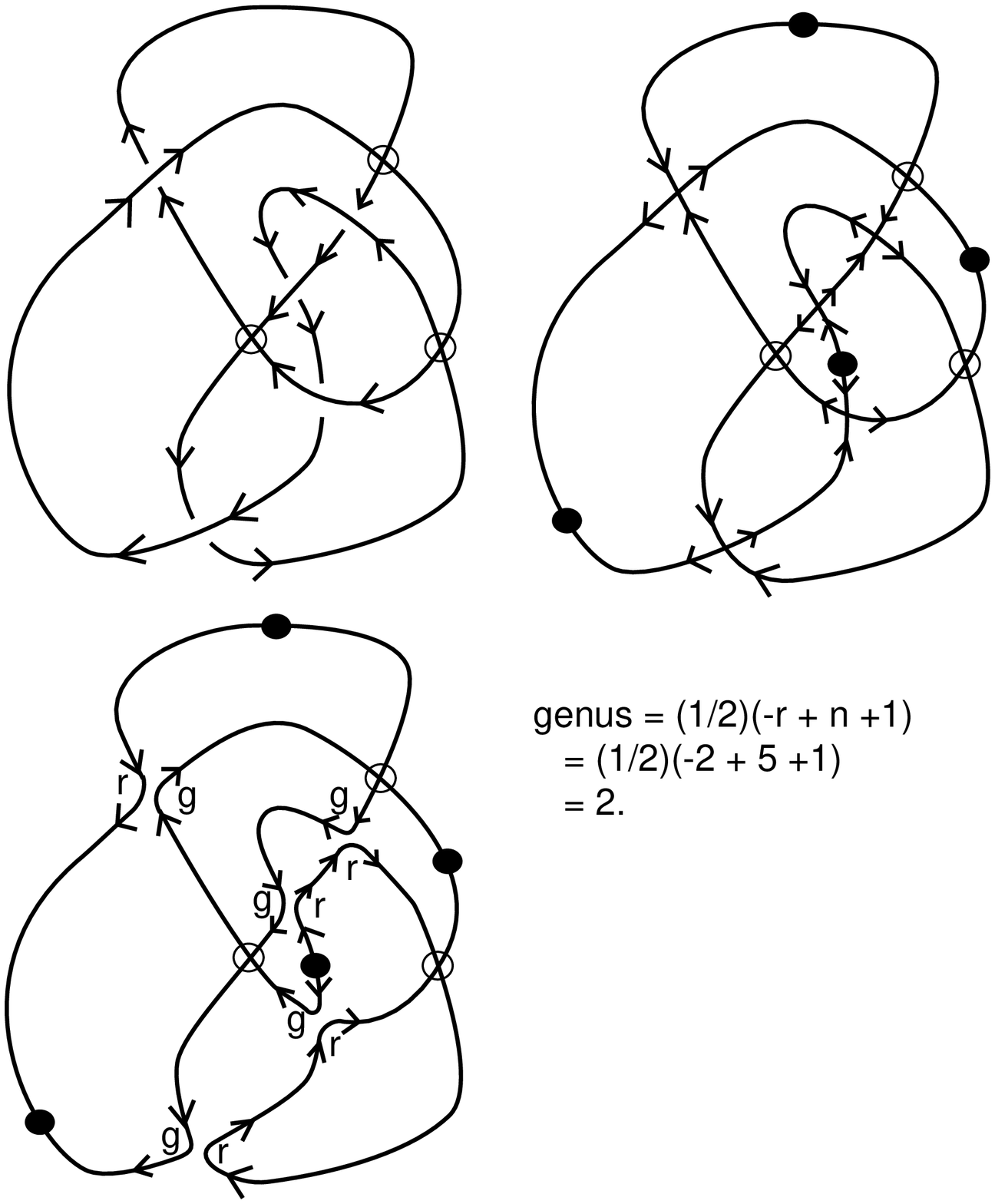}
     \end{tabular}
     \caption{\bf Lee Algebra Labels Seifert State for Specific Knot}
     \label{leeknot}
\end{center}
\end{figure}

\noindent{\bf Discussion.} This Theorem is proved by using a generalization of integral Khovanov homology to virtual knot theory devised by Manturov \cite{MBook}. In \cite{DKK} we reformulate this theory and show that it generalizes to the Lee homology theory (a variant of Khovanov homology) as well. In the Lee theory the chain complex is defined by a Frobenius algebra with the generators $g$ and $r$ as shown in Figure~\ref{leecut}. As this Figure and  Figure~\ref{leecycle} show, each link diagrm is reoriented in 
``source-sink" form so that every crossing has two arrows in and two arrows out that alternate cyclically around the crossing. In a virtual diagram it can be the case that these local orientations do not fit globally. There will be some edges where a switch of local orientation occurs from one end of the edge to the other. We mark points on those edges where a switch of orientation is required, and call these points the {\it cut locus}.  Loops in the state are labeled with algebra elements $g$ and $r$ so that $g$ is changed to $r$ and vice-versa when one move across the cut locus. Such labeled states are generators of the Lee Homology chain complex. It follows from the fact that $rg = gr = 0$ that such a labeled state represents an element of the kernel of the boundary mapping for the chain complex (hence a cycle in the complex) if an $r$ and a $g$ are paired at each smoothing site in the state.. Such labelings occur for Seifert smoothing states as is illustrated in
Figure~\ref{leecycle} and Figure~\ref{leeknot}. Since such cycles for positive knots occur with only 
$A$-type smoothings they cannot be in the image of the boundary map. Hence they represent the non-triviality of the positive (virtual) knot. Further analysis then yields the Theorem in direct analogy with Rasmussen's original arguments.
\bigbreak 

\bigbreak

\bigbreak

\bigbreak

\bigbreak

\subsection{Properties of the Virtual Stevedore's Knot}
We first point out that the virtual stevedore ($VS$) is an example that illustrates the viability of our theory.
We prove that $VS$ is not classical by showing that it is represented on a surface of genus one and no smaller. The reader should note the difference between representation of a virtual knot or link {\it on} a surface (as an embedding into the thickened surface) and the previous subsection's work on spanning surfaces.
\bigbreak

\begin{figure}
     \begin{center}
     \begin{tabular}{c}
     \includegraphics[width=7cm]{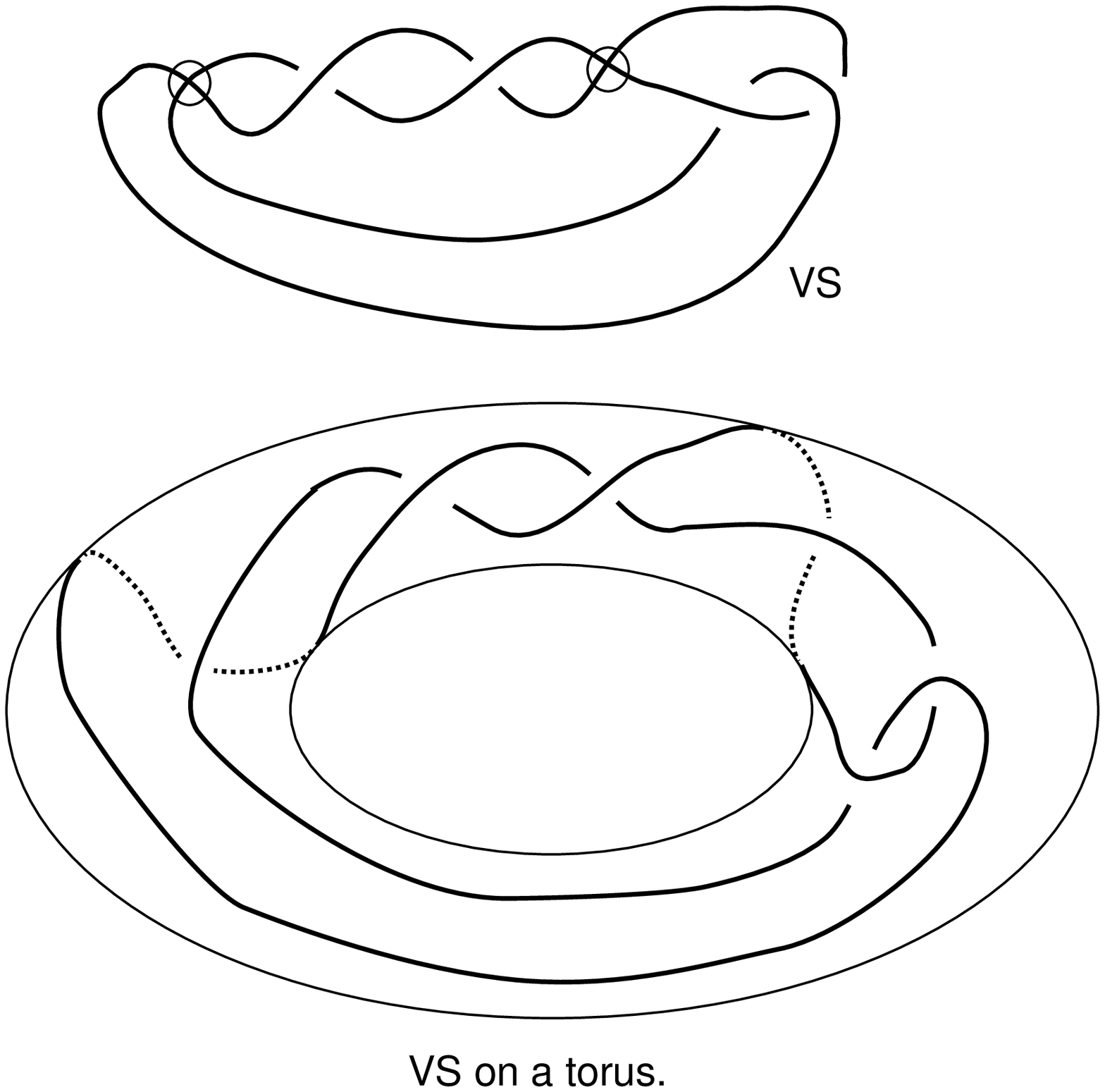}
     \end{tabular}
     \caption{\bf Virtual Stevedore on a Torus}
     \label{vsttorus}
\end{center}
\end{figure}

\begin{figure}
     \begin{center}
     \begin{tabular}{c}
     \includegraphics[width=7cm]{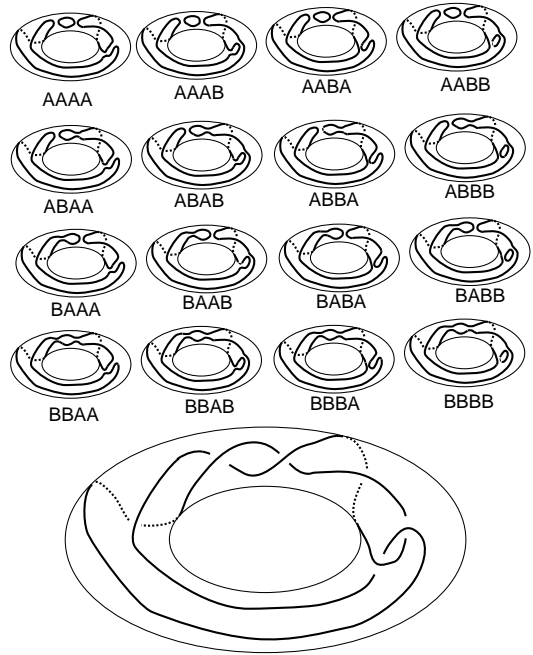}
     \end{tabular}
     \caption{\bf Virtual Stevedore is not Classical}
     \label{toral}
\end{center}
\end{figure}

The technique for finding this surface genus for the virtual stevedore  is to use the bracket expansion on a toral representative of $VS$ and examine the structure of the state loops on that surface. See Figure~\ref{vsttorus} and Figure~\ref{toral}. Note that in these Figures
the virtual crossings correspond to parts of the diagram that loop around the torus, and do not weave on the surface of the torus. An analysis of the homology classes of the state loops shows that the knot cannot be isotoped off the handle structure of the torus. See \cite{KaV4,MinSurf} for more information about using the surface bracket.
\bigbreak

\begin{figure}
     \begin{center}
     \begin{tabular}{c}
     \includegraphics[width=7cm]{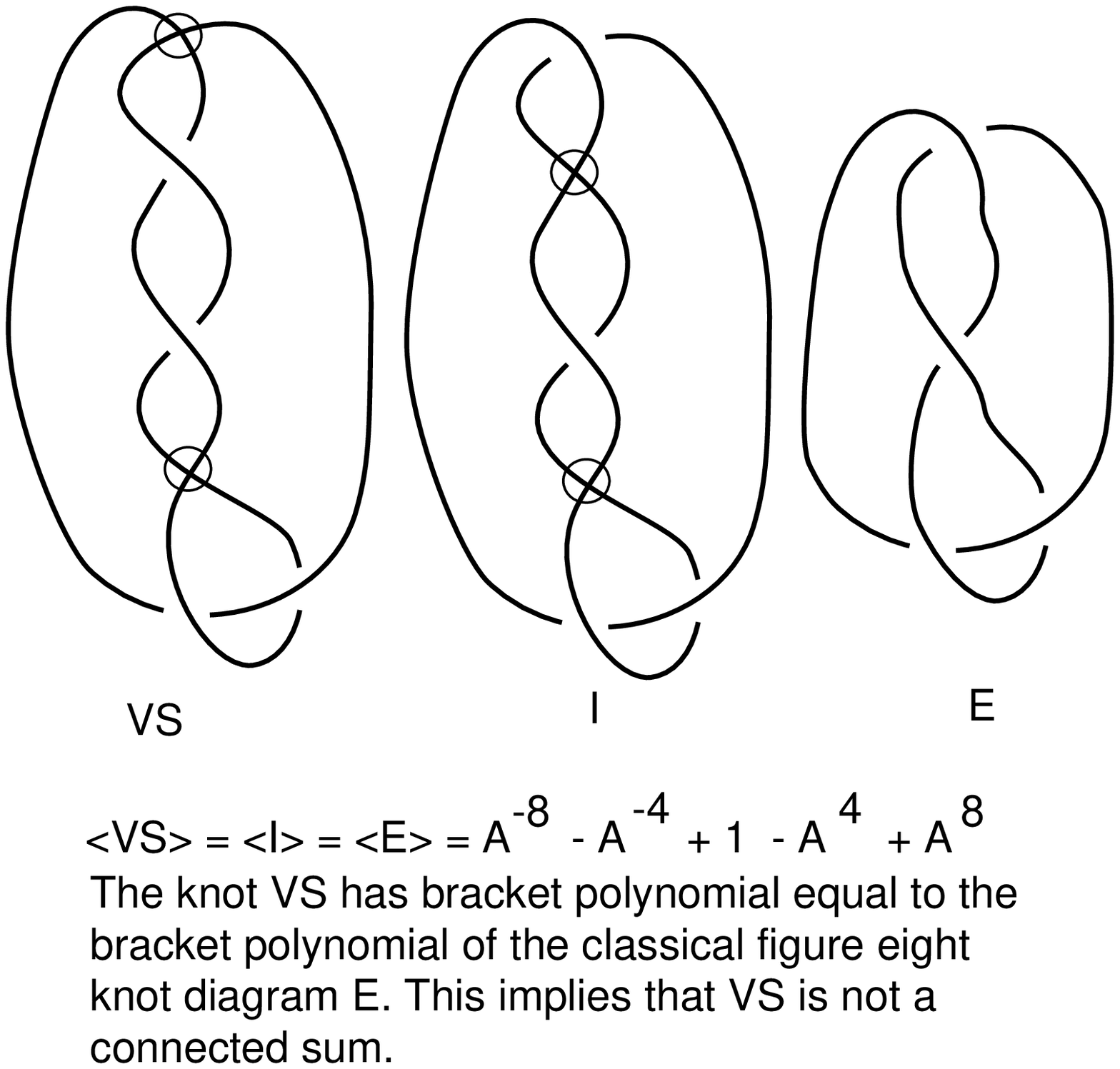}
     \end{tabular}
     \caption{\bf Bracket Polynomial of the Virtual Stevedore}
     \label{VSbracket}
\end{center}
\end{figure}

Next we examine the bracket polynomial of the virtual stevedore, and show as in Figure~\ref{VSbracket}
that it has the same bracket polynomial as the classical figure eight knot. The technique for showing this is to use the basic bracket identity for a crossing flanked by virtual crossings as discussed in the previous section. This calculation shows that $VS$ is not 
a connected sum of two virtual knots. Thus we know that $VS$ is a non-trivial example of a virtual slice knot. We now can state the problem:  {\bf Classify virtual knots up to concordance.} We will discuss this 
problem in this paper, but not solve it. The reader should note that the corresponding problem for classical knots is not solved, but has more techniques available. We will need to forge new techniques for the virtual problem.
\bigbreak

\subsection{Rotational Cobordisms}
In Section 1 we introduced rotational virtual knots, where the detour move is restricted to regular homotopy of plane curves. This means that the virtual curl of Figure~\ref{Figure 1} can not be directly simplified, but two opposite virtual curls can be created or destroyed by using the Whitney Trick of Figure~\ref{whitney}. We define {\it rotational virtual cobordism} by allowing saddles, deaths and births as before to extend the rotational equivalence relation. This section is a brief introduction to this theory in the form of two examples.
\bigbreak

\begin{figure}
     \begin{center}
     \begin{tabular}{c}
     \includegraphics[width=7cm]{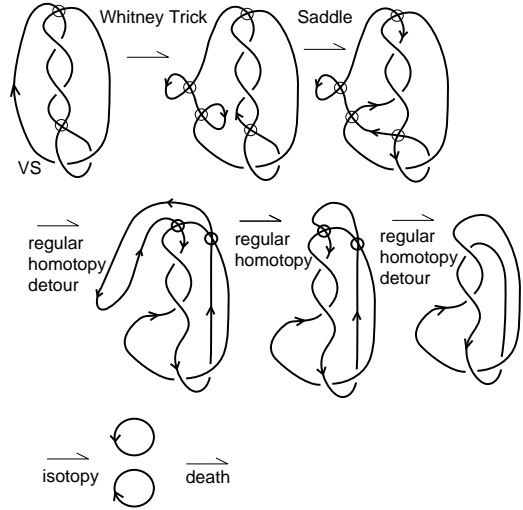}
     \end{tabular}
     \caption{\bf The Virtual Stevedore is Rotationally Slice}
     \label{rotslice}
\end{center}
\end{figure}

The first example is shown in Figure~\ref{rotslice}. We show that the virtual stevedore diagram is rotationally slice. This is accomplished by first doing a Whitney Trick, then using one curl, so produced, to produce a saddle. In this case there is a combination of regular homotopy and virtual isotopy (just labeled isotopy in the figure) that produces two disjoint circles with no virtual curls. These die to produce the virtual disk, showing that the knot is slice in the rotational category.
\bigbreak

\begin{figure}
     \begin{center}
     \begin{tabular}{c}
     \includegraphics[width=7cm]{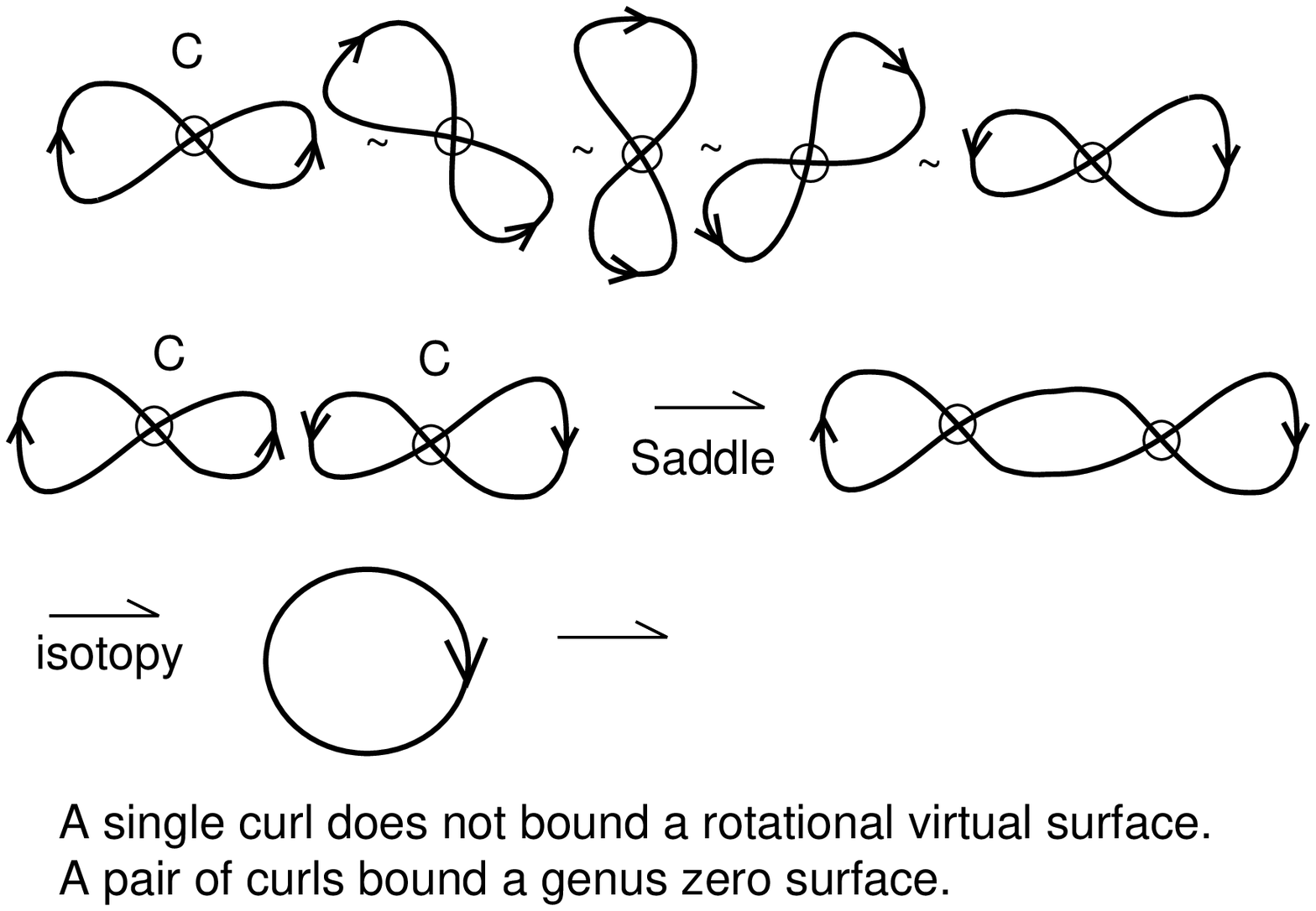}
     \end{tabular}
     \caption{\bf Parity and Cobordism of Curls}
     \label{curlcob}
\end{center}
\end{figure}

In Figure~\ref{curlcob} we illustrate how a single curl $C$ is (having Whitney degree zero) equivalent to its reverse (reversing orientation) and so can go through a saddle point with a copy of itself to a curve that is regularly homotopic to a trivial circle. Hence two curls can together rotationally bound a virtual surface of genus zero. This is a matter of parity! Given a rotational virtual link $L,$ define the 
{\it rotational parity of $L$}, $RotP(L)$ to be $0$ if $L$ has an even number of virtual crossings, and $1$ if $L$ has an odd number of virtual crossings. We have the 
\bigbreak

\noindent{\bf Theorem.} If $L$ is a rotational virtual link, then $L$ bounds a rotational virtual surface if 
and only if $L$ has even rotational parity, $RotP(L)= 0.$ 
\bigbreak

\noindent{\bf Proof.} Note that the parity of the number of virtual crossings is an invariant of rotational equivalence of virtual knots and links {\it and} it is also an invariant of saddle moves, death and births since these moves do not change the number of virtual crossings. We only allow deaths and births for circles that have no virtual crossings - call these {\it free circles}. Thus $RotP$ is an invariant of rotational cobordism. Therefore
if $RotP(L)$ is odd, $L$ cannot be cobordant to a disjoint union of free circles and so cannot rotationally bound a virtual surface. Conversely, if $RotP(L) = 0$ then the same cobordism we used before, a saddle move at every crossing, combined with isotopy, produces a collection of closed curves with only virtual crossings. These curves can be made into a disjoint collection of curves by regular isotopy, and then each curve is regularly isotopic to  curve with only ``external curls'', as illustrated in Figure~\ref{rotbound}. Then each such curve can undergo saddle moves to transform it to a disjoint union of curls of the form $C$ of Figure~\ref{curlcob}. This is also illustrated in Figure~\ref{rotbound}. The final number of curls of type $C$ is even since we assumed that the parity is even. Therefore the curls cancel in pairs as explained above, and we obtain a disjoint union of one-half their number as free circles. The circles bound disks. This finishes the construction of the surface and hence finishes this proof. //
\bigbreak

\begin{figure}
     \begin{center}
     \begin{tabular}{c}
     \includegraphics[width=7cm]{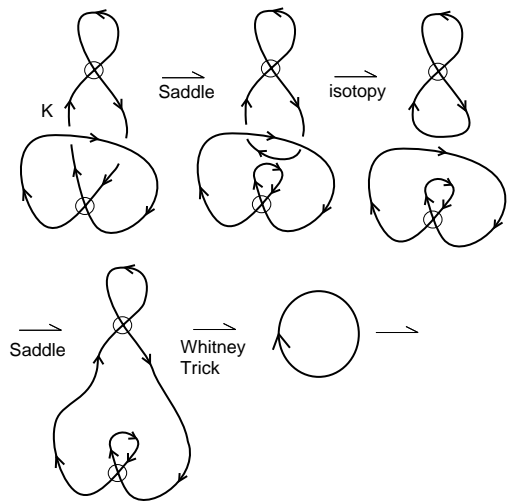}
     \end{tabular}
     \caption{\bf A Rotational Cobordism}
     \label{rotcob}
\end{center}
\end{figure}

\begin{figure}
     \begin{center}
     \begin{tabular}{c}
     \includegraphics[width=7cm]{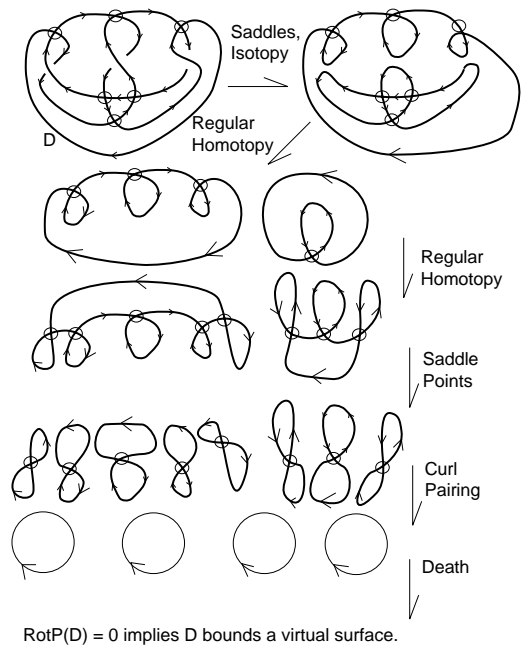}
     \end{tabular}
     \caption{\bf Diagram $D$ Bounds Rotational Virtual Surface if $RotP(D) = 0.$}
     \label{rotbound}
\end{center}
\end{figure}

We do not expect all virtual slice knots to be rotationally slice. Consider the rotational knot $K$ from Figure~\ref{rotateknot}. We proved in Section 1 that $K$ is a non-trivial non-classical rotational virtual knot. In Figure~\ref{rotcob} we illustrate a cobordism of $K$ to the disjoint union of two circles, each of which has  curl. Neither of these circles can bound a rotational virtual disc by our rules. So this cobordism stops short of exhibiting $K$ as a rotational slice knot. On the other hand, the two curls can interact through a saddle point to produce a free circle. Thus $K$ does rotationally bound a virtual surface of genus one. We conjecture that $K$ is not rotationally slice and that its least four-ball genus is one.
\bigbreak

Clearly much more work needs to be done in the study of cobordisms of rotational virtual knots and links.
Since there are many invariants of rotational knots and links (all the quantum link invariants), we can ask {\it how do quantum link invariants behave under rotational cobordism?} This will be the subject of subsequent papers.

\subsection{Band Passing}
\begin{figure}
     \begin{center}
     \begin{tabular}{c}
     \includegraphics[width=7cm]{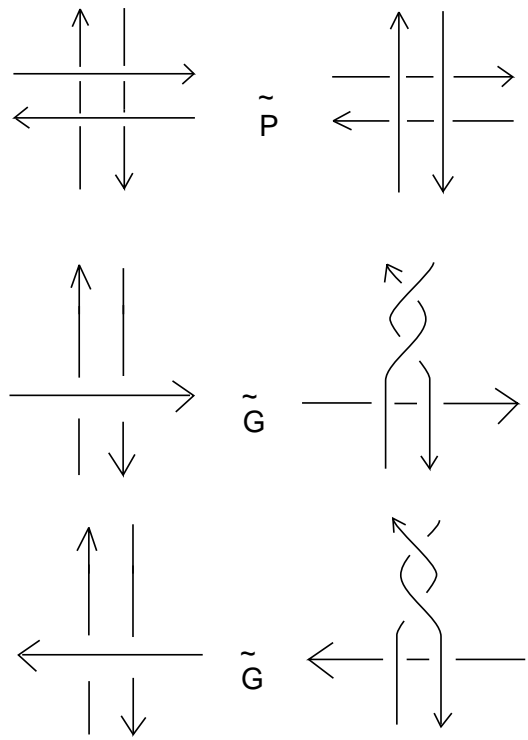}
     \end{tabular}
     \caption{\bf Pass and Gamma Moves}
     \label{passgamma}
\end{center}
\end{figure}

\begin{figure}
     \begin{center}
     \begin{tabular}{c}
     \includegraphics[width=7cm]{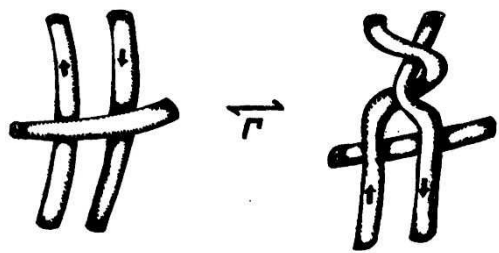}
     \end{tabular}
     \caption{\bf Gamma Move}
     \label{gamma}
\end{center}
\end{figure}

\begin{figure}
     \begin{center}
     \begin{tabular}{c}
     \includegraphics[width=7cm]{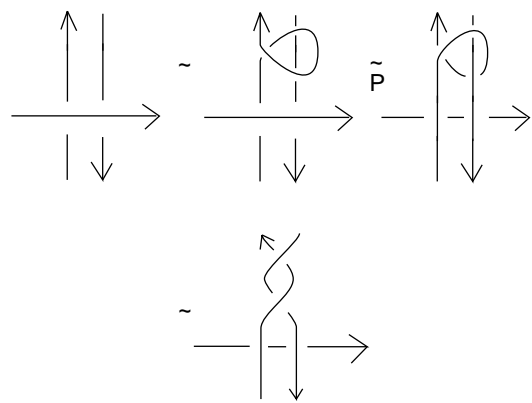}
     \end{tabular}
     \caption{\bf Gamma is Accomplished by Passing}
     \label{gammapass}
\end{center}
\end{figure}

\begin{figure}
     \begin{center}
     \begin{tabular}{c}
     \includegraphics[width=7cm]{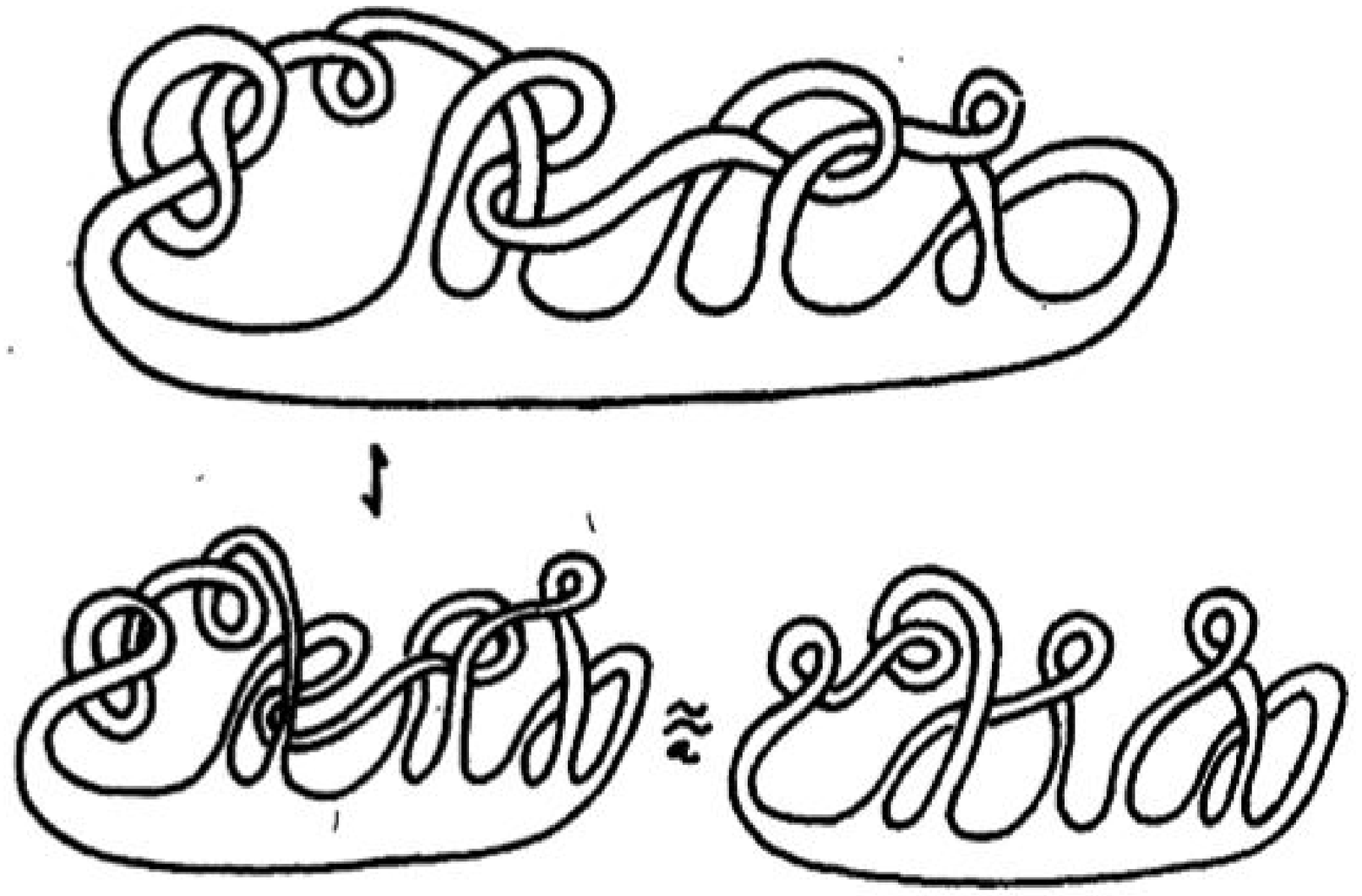}
     \end{tabular}
     \caption{\bf Simplifying a Surface by Passing Bands}
     \label{simplify}
\end{center}
\end{figure}

The Arf invariant of a  classical knot can be interpreted as the {\it pass-class} of the knot, where two knots are {\it pass-equivalent} \cite{OK} if one can be obtained from the other by ambient isotopy combined with switching pairs of oppositely oriented pairs of parallel strands as illustrated in Figure~\ref{passgamma}. The pass-class is a concordance invariant of classical knots and closely related to the Alexander polynomial. Any classical knot is pass-equivalent to either the trefoil knot or the unknot. The trefoil is pass-equivalent to its mirror image and is in a distinct pass-class from the unknot. The reader can get an idea of how this works for classical knots by examining Figure~\ref{simplify} where we show hows a complicated surface (with boundary a classical knot) can be simplified by band-passing. See 
\cite{OK} for more information about classical band passing.
\bigbreak

We would like to determine the pass-classes of virtual knots. This problem appears difficult at this time due the lack of invariants of the passing operation. We can obtain partial results by restricting passing to only odd crossings (crossings with an odd interstice in the Gauss code) but this is only a step on the way to understanding the pass equivalence relation for virtual knots. We expect that understanding this relation will shed light on problems of knot concordance.
\bigbreak

In Figure~\ref{passgamma} we illustrate pass-equivalence and also illustrate another move denoted by 
``G" in that figure and we refer to this move as the {\it gamma move.} The gamma move, illustrated separately in Figure~\ref{gamma},  switches one strand past two oppositely oriented strands and places a $2 \pi$ twist in these two strands. It is obvious that two gamma moves will accomplish a single pass-move, since the twist introduced by $G$-passing one strand is cancelled by the twist introduced by $G$-passing a second, oppositely oriented strand. It is also the case that any gamma move can be accomplished by a combination of ambient isotopy and a pass-move, as shown in 
Figure~\ref{gammapass}. Thus, pass-equivalence and gamma-equivalence are identical as equivalence relations on classical or virtual links. In the classical case, gamma-equivalence is of direct interest, as it is easy to see that a classical ribbon knot can be trivialized by gamma-moves (via the removal of ribbon singularities from an immersed disk that spans the link in three dimensional space). This is one way to see that classical ribbon knots are pass-equivalent to unknots. Thus, a classical knot that is not pass-equivalent to the unknot cannot be ribbon. 
\bigbreak

In the virtual case, we can see in some cases (such as the virtual stevedore's knot) that the knot is gamma-equivalent (hence pass-equivalent) to an unknot. But there exist examples of virtual knots that are slice, but are {\it not} pass-equivalent to the unknot. For example, consider the Kishino knot of 
Figure~\ref{kishino}. The Kishino has two oppositely oriented parallel strands in its middle, and one saddle point move transforms the Kishino into two virtual unknots. Hence the Kishino is slice. However,
we know (e.g. via the parity bracket) that the Kishino overlies a non-trivial flat virtual diagram. No pass-move changes the underlying flat diagram of a virtual knot. Hence pass-equivalence cannot unknot the Kishino diagram. This shows that there is a sharp difference between pass-equivalence for classical knots and pass-equivalence for virtual knots.
\bigbreak

\section{Virtual Surfaces in Four-Space}
We now define a theory of virtual surfaces in four-space that is given by moves on planar diagrams.
One of the projects of this proposal is to investigate the relationships between this diagrammatic 
definition and more geometric approaches to virtual $2$-knots due to Jonathan Schneider  and to 
Takeda \cite{Takeda}. We make diagrammatic definitions as follows: We use middle level markers as indicated
in Figure~\ref{markers} to encode two directions of smoothing a marked crossing in a planar diagram.
The classical interpretation of such a marker is that it represents a cobordism through a saddle point at 
the middle level ($t=0$ in the Figure) where the forms of smoothing above ($t=1$) and below ($t=-1$)
are shown via the conventions in the Figure. A diagram with markers can then be interpreted as two cobordisms attached at the middle. One cobordism goes downward to a collection of possibly linked and knottted loops, the other goes upward to another collection of linked and knotted loops. We will refer to these as the {\it up-cobordism} and the {\it down-cobordism}. A marked diagram is said to be 
{\it excellent} if both the up and the down cobordisms end in collections of unlinked circles that can be capped off with births (from the bottom) and deaths (at the top). The resulting schema is then a two-sphere and classically represents a two-sphere in four space. We take exactly this definition for a 
{\it virtual two-sphere} where it is understood that the ends of the two cobordisms will be trivial virtual links.
\bigbreak

\begin{figure}
     \begin{center}
     \begin{tabular}{c}
     \includegraphics[width=7cm]{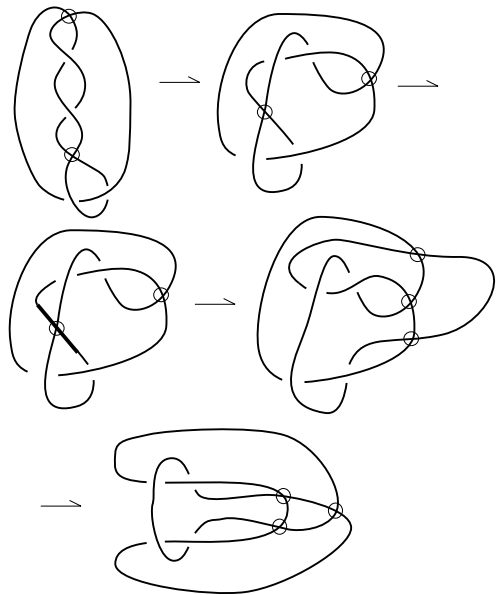}
     \end{tabular}
     \caption{\bf Converting Virtual Stevedore to a ``Ribbon Diagram"}
     \label{ribbondiagram}
\end{center}
\end{figure}

\begin{figure}
     \begin{center}
     \begin{tabular}{c}
     \includegraphics[width=6cm]{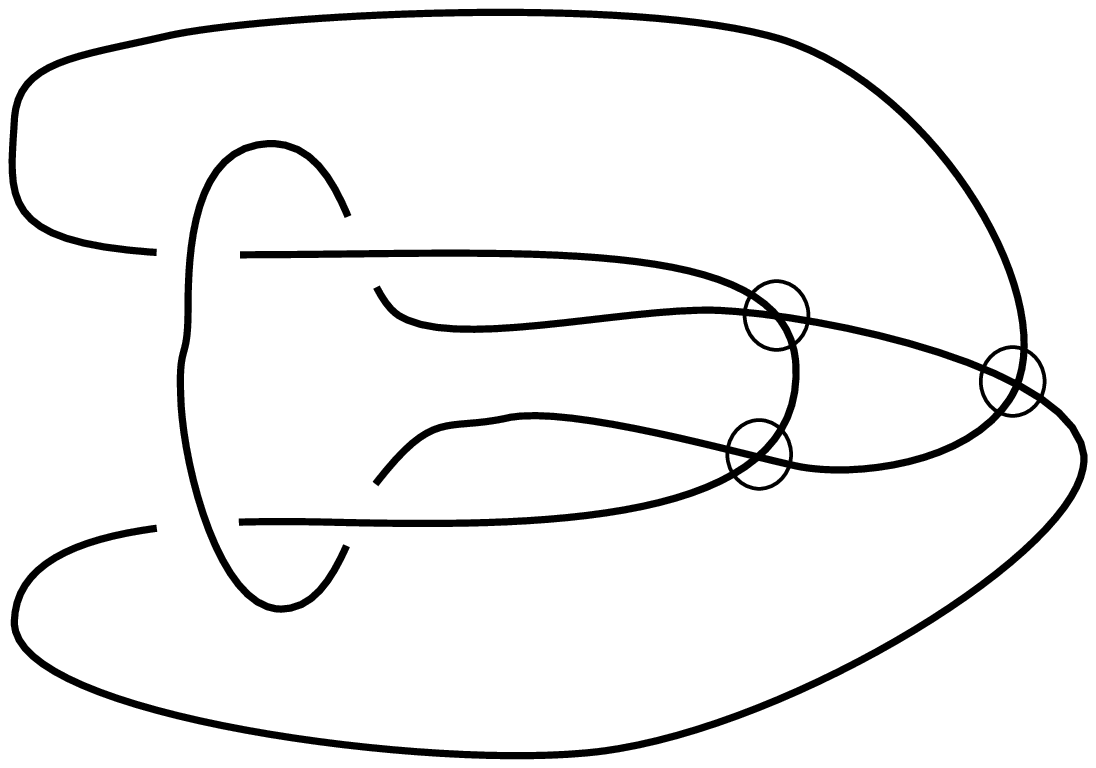}
     \end{tabular}
     \caption{\bf Ribbon Version of Virtual Stevedore}
     \label{ribbonstevedore}
\end{center}
\end{figure}

\begin{figure}
     \begin{center}
     \begin{tabular}{c}
     \includegraphics[width=6cm]{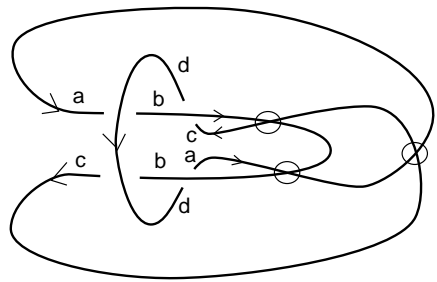}
     \end{tabular}
     \caption{\bf Labeled Ribbon Version of Virtual Stevedore}
     \label{labeledribbonstevedore}
\end{center}
\end{figure}

Just as in classical theory, if a virtual knot is slice, then we can make a virtual two-sphere from it by using the same cobordism both up and down. In Figure~\ref{dblslicesurface} we indicate the schema for such a surface involving one saddle point up and one saddle point down. Births for the original cobordism have to be represented directly in the middle level. The slicing example for $VS$, the virtual stevedore's knot, can be made into a two-sphere this way. We show the middle level diagram for this sphere, called $S$ in Figure~\ref{twotwo}.
In this same Figure, we show another middle level diagram for a virtual two-sphere S'. In this case we have used the fact (the reader can verify) that $VS$ can be sliced from its right-hand side. The sphere
S' is obtained by slicing upward from the left and downward from the right.
\bigbreak

We give moves on the middle level diagrams to define isotopy of the virtual two-spheres obtained from the middle level diagrams. The moves are indicated in Figure~\ref{middlemoves}. They are a virtual
generalization of the Yoshikawa moves that have been studied\cite{Swenton,SYLee}  for isotopies of the classical middle level formulations. Thus we say the two two-spheres are {\it isotopic} if one can be obtained from the other via these {\it Generalized Yoshikawa moves}. In particular, the fundamental group of the two-sphere,defined by adding relations at saddle points exactly as in the classical case (but from the virtual knot theoretic fundamental group) is an isotopy invariant. For example, in Figure~\ref{fundgrp} we calculate the fundamental group of $VS$ and find that, in it the arcs whose elements must be identified to obtain the fundamental group of the sphere $S$ of Figure~\ref{twotwo} are already identified in the fundamental group of $VS$. Thus we find that the sphere $S$ is knotted since it has the same non-trivial fundamental group as $VS.$ On the other hand, it is not hard to see that the fundamental group of the sphere $S$ is isomorphic to the integers. At this writing I do not know if this sphere is virtually unknotted.

\begin{figure}
     \begin{center}
     \begin{tabular}{c}
     \includegraphics[width=6cm]{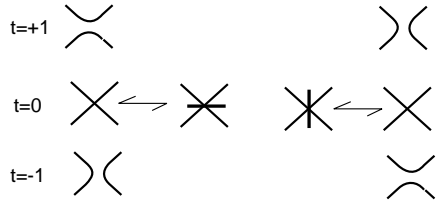}
     \end{tabular}
     \caption{\bf Middle Level Markers}
     \label{markers}
\end{center}
\end{figure}

\begin{figure}
     \begin{center}
     \begin{tabular}{c}
     \includegraphics[width=7cm]{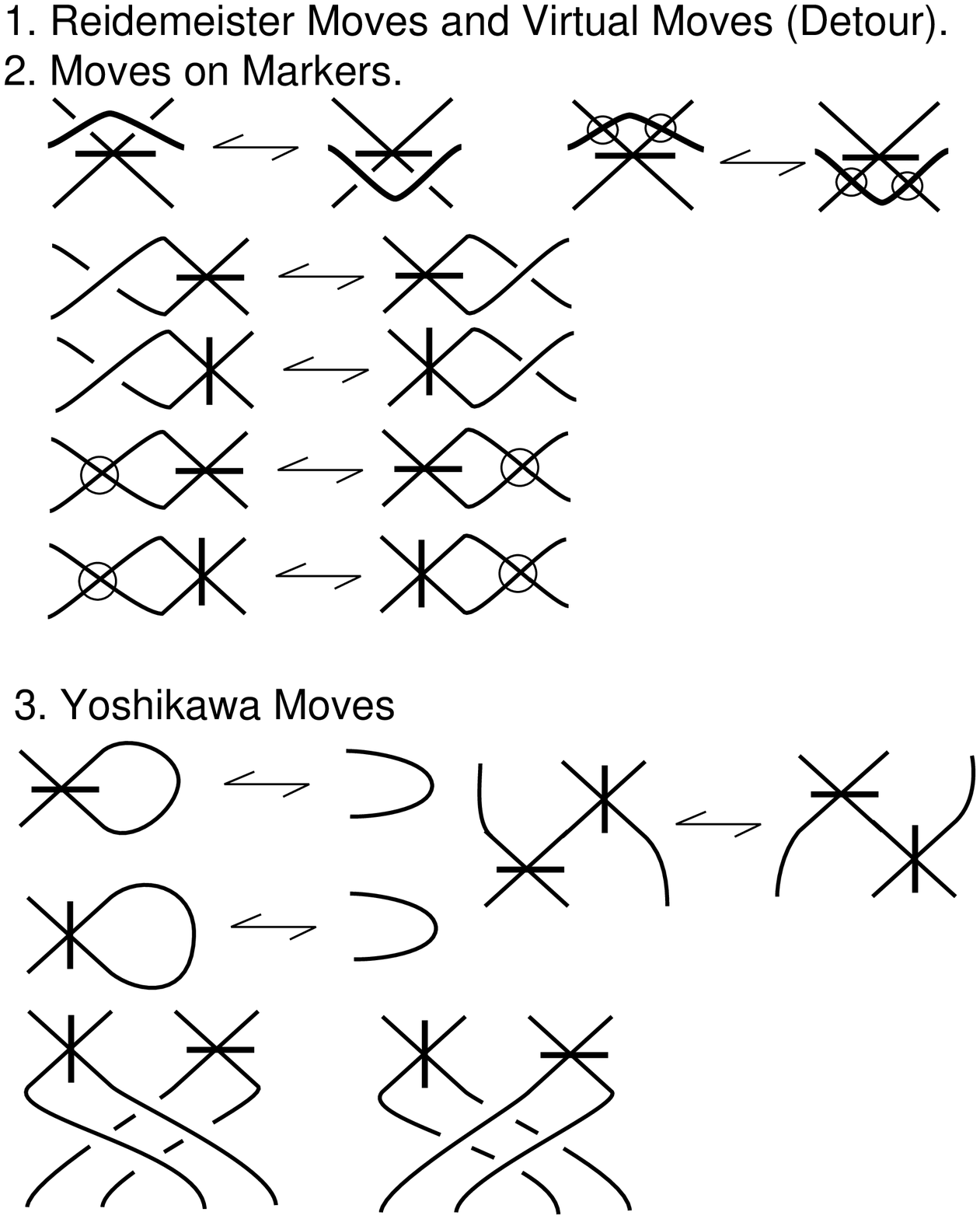}
     \end{tabular}
     \caption{\bf Middle Level Moves}
     \label{middlemoves}
\end{center}
\end{figure}

\begin{figure}
     \begin{center}
     \begin{tabular}{c}
     \includegraphics[width=7cm]{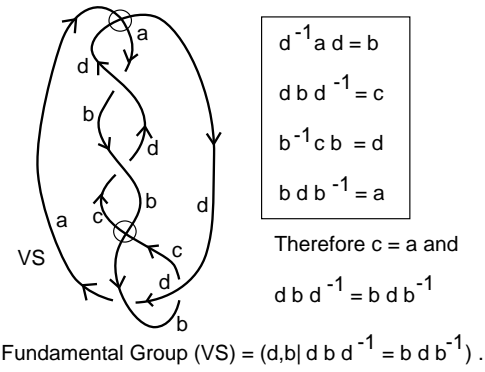}
     \end{tabular}
     \caption{\bf Fundamental Group of $VS$}
     \label{fundgrp}
\end{center}
\end{figure}

\begin{figure}
     \begin{center}
     \begin{tabular}{c}
     \includegraphics[width=3cm]{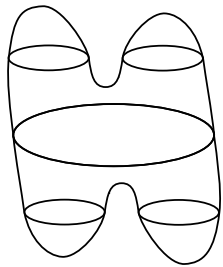}
     \end{tabular}
     \caption{\bf Abstract Double Slice Surface}
     \label{dblslicesurface}
\end{center}
\end{figure}

\begin{figure}
     \begin{center}
     \begin{tabular}{c}
     \includegraphics[width=7cm]{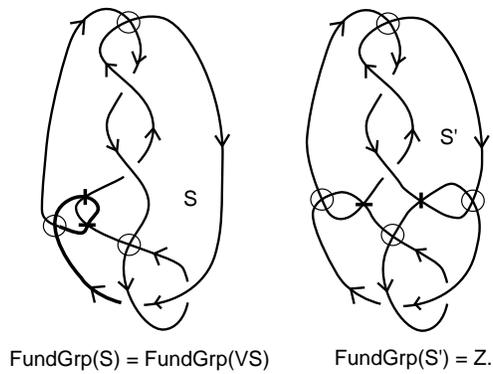}
     \end{tabular}
     \caption{\bf Two Two-Spheres}
     \label{twotwo}
\end{center}
\end{figure}

The generalized Yoshikawa moves present a useful first formulation for a theory of virtual surfaces.
One of the advantages of this approach is that we can adapt the generalization of the bracket polynomial of Sang Youl Lee \cite{SYLee} to obtain a bracket invariant for virtual two-spheres. This will be an important  subject of investigation for this proposal. We want to know how this diagrammatic formulation is related to immersions of surfaces in four space that could represent virtual two-knots. In this case the levels 
(movie of a cobordism) description that we have adopted gives such an immersion, and one can begin the investigation at that point. For these reasons, we believe that this formulation of virtual cobordism an virtual surfaces will be very fruitful and lead to many new results.

\bigbreak

\end{document}